\documentclass[preprint,12pt,authoryear]{elsarticle}
\usepackage{graphicx,lmodern,amsfonts,mathrsfs,amsmath,amssymb,latexsym,xspace,epsfig,psfrag,color,tikz,tikzscale,pgfplots,colortbl}
\usepackage[colorlinks,bookmarksnumbered,citecolor=blue]{hyperref}
\usepackage[ruled]{algorithm2e}
\usepackage[utf8]{inputenc}  
\usepackage{scalefnt}
\usepackage{subfigure}
\usepackage{multirow}
\usepackage[flushleft]{threeparttable}
\usepackage{pdflscape}
\usepackage{multicol}
\usepackage{natbib}
\usepackage{hyperref}
\usepackage{lineno}
\usepackage{longtable, lscape} 
\usepackage{bbm}
\usepackage{calrsfs}
\usepackage{txfonts}
\usepackage{rotating}
\usepackage{comment}

\usepackage[top=30mm,bottom=20mm,left=30mm,right=20mm,pdftex,includeheadfoot]{geometry}

\usepackage{booktabs}
\usepackage{multirow}
\usepackage{pgfplotstable}
\usepackage{filecontents}
\usepackage{xspace}
\usepackage{tabularx}
\usepackage{siunitx}

\usepackage{tikz,pgfplots}
\pgfplotsset{compat=1.15}
\usepgfplotslibrary{statistics}
\usetikzlibrary{shapes,arrows}

\definecolor{Gray}{gray}{0.95}
\definecolor{Gray1}{gray}{0.7}
\definecolor{Gray2}{gray}{0.8}
\definecolor{Gray3}{gray}{0.9}
\definecolor{Gray4}{gray}{1}

\journal{European Journal of Operational Research}

\begin{document}

\begin{frontmatter}

\title{A hybrid adaptive Iterated Local Search with diversification control to the Capacitated Vehicle Routing Problem}

\author{Vinícius R. Máximo}
\ead{vinymax10@gmail.com}

\author{Mariá C. V. Nascimento\corref{cor1}}
\ead{mcv.nascimento@unifesp.br}
\cortext[cor1]{Corresponding author}

\address{Instituto de Ciência e Tecnologia, Universidade Federal de São Paulo (UNIFESP)\\ Av. Cesare M. G. Lattes, 1201, Eugênio de Mello, São José dos Campos-SP, CEP: 12247-014, Brasil
}

\begin{abstract}
Metaheuristics are widely employed to solve hard optimization problems, like vehicle routing problems (VRP), for which exact solution methods are impractical. In particular, local search-based metaheuristics have been successfully applied to the capacitated VRP (CVRP). The CVRP aims at defining the minimum-cost delivery routes for a given set of identical vehicles since each vehicle only travels one route and there is a single (central) depot. The best metaheuristics to the CVRP avoid getting stuck in local optima by embedding specific hill-climbing mechanisms such as diversification strategies into the solution methods. This paper introduces a hybridization of a novel adaptive version of Iterated Local Search with Path-Relinking (AILS-PR) to the CVRP. The major contribution of this paper is an automatic mechanism to control the diversity step of the metaheuristic to allow it to escape from local optima. The results of experiments with 100 benchmark CVPR instances show that AILS-PR outperformed the state-of-the-art CVRP metaheuristics.
\end{abstract}


\end{frontmatter}


\section{Introduction}


The Capacitated Vehicle Routing Problem (CVRP) is a variant of the vehicle routing problem (VRP) and was initially proposed by \cite{Dantzig1959}. The CVRP consists of finding a set of routes that minimizes the costs of a fleet of homogenous vehicles in order to serve all customers respecting the capacity limit of the vehicles. There is a single depot from which all vehicles start and end their routes. The costs involved in the CVRP are symmetric and some authors explicit this feature by defining the problem as Symmetric CVRP (SCVRP). Problems that consider asymmetric costs are known as Asymmetric CVRP (ACVRP) \citep{Eksioglu2009}.

The pioneering exact methods for the CVRP were proposed by \cite{Laporte1987}. Some surveys on exact algorithms can be found in \citep{Toth1998,Toth2002,Baldacci2007,Cordeau2007,Baldacci2010,Baldacci2011,Poggi2014}. The most competitive  exact algorithms are based on  branch-and-cut-and-price  \citep{Fukasawa2006,Pecin2014}, being able to solve instances to up to 360 customers. The exact solution of the CVRP, however, requires a significantly large computational effort, the reason  why heuristic methods stand out from the solution methods \citep{Christiaens2020}. The Iterated Local Search (ILS) with Set Partitioning, called ILS-SP, proposed by \cite{Subramanian2013}, the Unified Hybrid Genetic Search (UHGS) proposed by \cite{Vidal2012, Vidal2014}  and the Slack Induction by String Removals (SISRs) proposed by \citet{Christiaens2020} are currently the state-of-the-art heuristic methods to the CVRP.  

\cite{Uchoa2017} recently proposed a new set of benchmark instances for the CVRP and compared the performance of ILS-SP and UHGS by carrying out experiments on these instances. The authors demonstrated that for instances with up to 200 vertices the algorithms perform reasonably well. However, for instances with more than 200 vertices, the problem remains a major computational challenge. SISRs were also tested with this new set of benchmark instances and achieved better results than ILS-SP and UHGS for the largest instances  -- with 500 to 1000 vertices. 

ILS is a local-search-based metaheuristic applied successfully in a wide variety of combinatorial optimization problems. Broadly, ILS consists of two phases: the perturbation phase and the local search phase. The perturbation phase is applied to a high-quality local optimum, the so-called reference solution. It aims at obtaining a solution outside the neighborhood of the reference solution in order to explore new regions of the search space. The local search phase is an improvement strategy applied to the solution while a local optimum is not found. These two steps are repeated until a stop criterion is reached. The degree of perturbation of solutions in ILS is of paramount importance for the performance of the algorithm \citep{Lourenco2003}. The main problem to be solved by the diversity control mechanism is the identification of the optimal perturbation intensity to be applied to the solution for the algorithm to perform better. Another important aspect approached in ILS is the solution's acceptance criterion of the method. The success of hill-climbing strategies highly depends on the decision regarding which solutions will have their neighborhoods investigated.

In line with this, in this paper, we present a hybridization of an adaptive version of ILS \citep{Lourenco2003,Lourencco2019} with Path-ReLinking (PR) \citep{Glover1996}  to the CVRP.  The major contribution of this study is the diversity control mechanism in the method, by guiding the perturbation degree and the acceptance criterion by  adaptive strategies.  The proposed strategy to address the perturbation intensity promotes this indication at execution time. For this, its value is guided by the distances between the reference solutions and the solutions found after the local search, obtained in earlier iterations of the method.  In addition, we present a criterion that takes into account the flow of accepted solutions during the search process to meet the acceptance criterion. A threshold that restricts the solutions whose neighborhood will be further investigated is estimated considering the flow evaluated in execution time.
Computational experiments were carried out with the benchmark instances proposed in \cite{Uchoa2017}. The adaptive ILS PR-hybrid, here named AILS-PR, achieved superior performance in relation to the state-of-the-art heuristic methods for the CVRP: UHGS, ILS-SP and SISRs. In comparison to the three algorithms,  AILS-PR presented the lowest mean gaps in 93\% of the tested  instances.

The remainder of this paper is organized as follows. Section~\ref{sec:problem} presents the problem description and notation used throughout the paper. Section~\ref{sec:related} briefly reviews studies closely related to this paper. Section \ref{sec:proposedAILS} introduces a new adaptive ILS, named AILS, to solve optimization problems. In addition, this section presents the proposed AILS to solve the CVRP. Section \ref{section:ILSPR} shows a hybrid AILS with PR to approach the CVRP. Computational experiments are shown in Section \ref{section:EC}, as well as the analysis of the performance of AILS-PR in comparison to the state-of-the-art heuristic solution methods for the CVRP. Section \ref{section:conclusao} sums up the main contributions of this paper and presents the final remarks.

\section{Problem description}\label{sec:problem}


The CVRP is presented in the following form: Let $G=(V,E)$ be an undirected complete graph  where $ V = \{0,1, \ldots, n\}$ is the set of $ n + 1 $ vertices and $ E $ the set of edges. The depot is represented by  vertex $ 0 $ and $V_c = V\setminus \{0\}$  represents the $n$ clients to be visited. Each edge $( i,j) \in E=\{(i,j):i,j\in V, i<j\}$ is associated with a non-negative weight $d_{i,j}$ which indicates the cost related to moving from point $i$ to point $j$. Such costs are symmetric, that is, $d_{i,j} = d_{j,i} ~\forall i,j \in V$.

In practical cases, the weights on edges may represent the distance between two points, the travel time, or the cost of the trip. Each customer $ i \in V_c $ has a non-negative demand $ q_i $ that must be met, and the deposit has demand $ q_0 = 0 $. To meet customers' demands, the fleet is composed of $ m $ identical vehicles limited to carry $ \bar{q} $ units of demands. The demand of each customer is smaller than the vehicle's capacity, that is, $q_i \leq \bar{q}, ~\forall i \in V_c$. 
 
The CVRP consists of finding a set of $ m$ closed walks that start and end at the depot. Moreover, the closed walks must minimize the sum of the costs of each closed walk and respect the following constraints:

\begin{itemize}
\item The sum of demands of the nodes of each closed walk does not exceed the vehicle capacity. 
\item Each vertex must belong to a single closed walk.
\end{itemize}

Each route is represented by a closed walk without node repetition.  A closed walk is represented by a cyclic sequence of vertices where a pair of vertices is adjacent if they are consecutive in the sequence and non-adjacent, otherwise. The routes of a solution $s$ are described by $\mathscr{R}=\{R^s_1,R^s_2,\ldots, R^s_m\}$, where $R^s_i = \{v^{i}_0, v^{i}_1,\ldots,v^{i}_{m_i}\}$, where $m_i$ is the size of route $R^s_i$,  $v^{i}_0=v^{i}_{m_i}=0$ and $R^s_i \cap R^s_j =\{0\}$, when $i\neq j$ and $\cup_{i=1}^m R^s_i  = V$.

The CVRP belongs to the class of  $\mathbf{NP}$-hard problems \citep{Lenstra1981}. Heuristic methods are widely employed to solve challenging instances even though the existence of exact methods. In the next section, we present a brief literature review on the state-of-the-art methods for the CVRP.

\section{Related Works}\label{sec:related}

This section presents the related  literature  to the proposed study. First, the current state-of-the-art metaheuristics for the CVRP are briefly showed. Then, an overview of adaptive Iterated Local Search is introduced.

\subsection{Metaheuristics for the CVRP}

The literature on heuristics and meta-heuristics for the CVRP is substantial. In particular, studies on classical heuristics and metaheuristics up to 2014 were summarized in the surveys presented by \cite{Laporte2002,Gendreau2001,Gendreau2008, Laporte2014}. A recent review on metaheuristic algorithms   proposed in 2009 to 2017 to solve the VRP is presented by \cite{Elshaer2020}. The review shows a taxonomy of the metaheuristics to the VRP and its variants. According to the authors,  98.9\% of the 299 reviewed papers tackle a CVRP variant. Because of the extensive literature on metaheuristics for the CVRP, in this section, we only focus on the current state-of-the-art heuristic methods.

In \citep{Subramanian2013}, the authors introduced ILS-SP, a hybridization of the ILS metaheuristic with an exact algorithm that considers the set partitioning (SP) formulation to solve the CVRP. In this solution method, the best routes generated by ILS throughout the search process are stored. Then, the exact algorithm attempts to solve the partitioning problem by choosing the optimal set of non-overlapping routes whose union is the whole set of vertices.  Successive swap and shift movements and a  Randomized Variable Neighborhood Descent  (RVND) are the solution perturbation strategies  in ILS. In RVND, the best improvement strategy is applied to each neighborhood -- 2-opt and 2-swap -- until no further improvement is observed.  ILS-SP is able to solve several variants of VRP, according to the authors. The results of computational experiments performed with hundreds of instances of the VRP variants showed that ILS-SP performed well in instances with up to 480 vertices.

\cite{Vidal2014} presented a heuristic method named  Unified Hybrid Genetic Search (UHGS) which solved 29 CVRP variants. UHGS has a local search phase that employs a 2-exchange and 2-opt neighborhoods. The individuals of the genetic algorithm are represented without route separation and use the same optimal route definition algorithm presented by \citet{Prins2004}. This representation allows the use of simpler crossing operators, the same suggested by \cite{Prins2004}, being the selection operator a binary roulette wheel. According to the authors, allowing unfeasible solutions contributes notably to the good performance of the method. This algorithm also uses a diversity control, thus avoiding premature convergence of the method. In this way, each individual is evaluated according to the solution cost penalized by the violated constraints and the distance in relation to the other individuals in the population. Considering CVRP, UHGS was tested in a new set of benchmark instances proposed by \citet{Uchoa2017}. It obtained good results in comparison to the ILS-SP proposed by \citet{Subramanian2013}, being state-of-the-art for CVRP. 

\citet{Christiaens2020} proposed a method   called Slack Induction by String Removals (SISRs), which has the same methodology as the Adaptive Large Neighborhood Search (ALNS). The authors chose a simple version of the ruin-and-recreate (R\&R) algorithm, using only one ruin method and one recreate procedure. The ruin method consists of removing some sequences of adjacent vertices belonging to the same route. The recreate method inserts the removed vertices in the lowest cost position, however, some positions may not be verified. The authors call this non-verification blinks that occur with a certain probability. The performance of SISRs was compared with UHGS \citep{Vidal2014} and ILS-SP \citep{Subramanian2013} using the 100 instances introduced in \cite{Uchoa2017}.  SISRs presented a better average solution in 57 instances. Besides, the average time was lower in comparison to the other algorithms. 

\subsection{Adaptive versions of the Iterative Local Search}

Iterated Local Search (ILS) is an iterative local search algorithm proposed by \citep{Lourenco2003} to solve combinatorial optimization problems. In the already mentioned review on metaheuristics for the VRP, \cite{Elshaer2020} state that about 10\% of the metaheuristics proposed in the reviewed literature to the VRP is based on the ILS. In particular, to the CVRP, the authors cite the works  \citep{Chen2010,Cordeau2012,Subramanian2013}. In the three studies, the authors either hybridize ILS with another heuristic or exact method to a broad exploration of the search space. 

As discussed earlier, the perturbation degree to escape from local optima is an important step of the method. In the traditional ILS, the perturbation intensity is randomly chosen at each iteration of the method. However, to some applications, to take such a decision randomly may not be the best choice. In line with this, a few studies suggest adaptive strategies to adjust the perturbation intensity in the literature, as, for example in \citep{Dong2015}.

\cite{Dong2015} proposed an adaptive version of  ILS to approach the flow shop scheduling problem. In their solution method, a mechanism to identify the appropriate perturbation strength basing on the status of the neighbor solutions found in the local search is proposed. If the average quality of the neighbor solutions decreases, then the strategy intensifies the perturbation strength to force the current solution to move to another region of the search space.

Another adaptive systematic in ILS  found in the literature is based on the dynamic prediction of the quality improvement based on past iterations for choosing the perturbation heuristics. In line with this, \cite{Walker2012}  suggested an adaptive ILS to the CVRP with time windows. The strategy proposed by the  authors employs online learning mechanisms to select the movements of the perturbation step of the algorithm. \cite{Schneider2015} proposed an adaptive VNS algorithm for VRPs with intermediate stops that chooses the perturbation heuristics basing on probabilities periodically adjusted along the search process.  To solve the CVRP, we tested a similar probabilistic approach, which estimates the probabilities to choose the perturbation heuristics based on the quality improvement information collected in past iterations.   
However, we opted for randomly choosing  which perturbation heuristic to apply during the search process since the probabilistic-based approach did not present better results.

\section{Adaptive Iterated Local Search (AILS)}\label{sec:proposedAILS}

In the proposed Adaptive Iterated Local Search, referred to as AILS, the perturbation degree and the acceptance criterion rely on adaptive strategies. In summary, the  perturbation degree intensity is iteratively adjusted bearing on the distance between the last solution found after the local search and a reference solution. The adaptive strategy to define the acceptation criterion is based on a threshold value to avoid the acceptance of low-quality solutions,  dynamically adjusted along the iterations. The primary idea behind the adjustment lies in the empirical observation that promising regions are usually composed of numerous good quality solutions. In this case, the algorithm must curb the acceptance criterion, for the algorithm to intensify the search in that region. The number of good solutions in non-promising regions is usually low and the algorithm must relax the acceptance criterion for the algorithm to escape from such a region.

Algorithm~\ref{AILS} shows a general framework of the introduced Adaptive Iterated Local Search (AILS).

\begin{algorithm}[!htb]
\LinesNumbered
\SetAlgoLined
\KwData{Instance data}
\KwResult{The best solution found $s^*$}
$s \leftarrow$ Construct an initial solution \\ 
$s^r,s^* \leftarrow$ Local Search($s$) \\
\Repeat{stop criterion is met}
{ 
    $s \leftarrow$ Perturbation Procedure  ($s^r$, $\mathscr{H}^r_k$)\\
    $s \leftarrow$ Local Search($s$)\label{linha5ails}\\
	Update the diversity control parameter $\omega_{\mathscr{H}^r_k}$ considering the distance between $s$ and $s^r$\\
    $s^r\leftarrow $ Apply acceptation criterion to $s$\\
	Update the acceptance criterion\\
	Assign $s$ to $s^*$ if $ f(s) < f(s^*)$\\

}
\caption{Adaptive Iterative Local Search}
\label{AILS}
\end{algorithm}

In Algorithm~\ref{AILS}, an initial solution $s$ is constructed for the studied problem.  A local search is applied to $s$ and the resulting local optimum is assigned to $s^*$ and $s^r$, respectively, the best overall solution and the reference solution.  The iterative  process begins by first applying to the reference solution a perturbation heuristic $\mathscr{H}^r_k$, where $k\in\{1,\ldots,n_k\}$, being $n_k$ the number of perturbation heuristics of the strategy. As told before, our approach employs the random selection of $k$ at each iteration of AILS. Then, the solution resulting from the perturbation is called $s$, to which the local search is applied. After that, as part of the proposed adaptive strategy,   the value of the diversity control parameter $\omega_{\mathscr{H}^r_k}$, which is dependent on the perturbation heuristic $\mathscr{H}^r_k$, is updated -- thoroughly explained in the next section. Thus,  the acceptation criterion decides if solution $s$  replaces the current reference solution $s^r$. Again, as part of the adaptive framework here proposed, the acceptation criterion is updated as discussed in Section~\ref{criterioAceitacao}. The best overall solution is updated in the sequence. The iterative process stops when the stop criterion is met.

It is worth pointing out that this algorithm has  the main structure of the classical ILS, except for lines 6 and 8. In this paper, we introduce a diversity control strategy to define the perturbation strength as well as the acceptation criterion, both adaptive. Each of them is  explained at length in the next sections.

\subsection{Diversity control}
\label{ControleDiversidade}

The diversity control is highly important for the algorithm to escape from a local optimum region. The proposed algorithm has two mechanisms to control diversity in the search process. The first is related to the degree of perturbation of the current solution. The second controls the diversity through the acceptance criterion. 

\subsubsection{Perturbation Control Method}
\label{ControleGrauPertubacao}

The perturbation degree control has the objective of adjusting the intensity of the solution perturbation in ILS regardless of the instance being solved.

To indicate the degree or extent of a given perturbation, it is necessary to measure the distance between the reference solution $s^r$ and the solution $s$ obtained after the local search. For example, to approach the CVRP, we employ  the symmetric distance presented in Equation \eqref{dab}. In this case, let $E^s$ be the set of edges of a solution $s$. 

\begin{align}
& d(s,s^r) = |E^{s} \triangle E^{s^r}|& \label{dab}
\end{align}

The value of $\omega_{\mathscr{H}^r_k}$ determines the degree of perturbation to be employed by a perturbation heuristic $\mathscr{H}^r_k$. The higher its value, the greater the number of vertices that will be removed in the perturbation control method. Therefore, we set the value of $\omega_{\mathscr{H}^r_k}$ to ensure a controlled diversity in the search. For this, the adjustment process of $\omega_{\mathscr{H}^r_k}$ occurs as described in Algorithm \ref{updateOmega}.

\begin{algorithm}[!htb]
\LinesNumbered
\SetAlgoLined
\KwData{Solutions $s$ and $s^r$, $it_{\mathscr{H}^r_k}$ $d_{\mathscr{H}^r_k}$, $\gamma$, $\omega_{\mathscr{H}^r_k}$, $d_{\beta}$ and $size$.}
\KwResult{ Updated $\omega_{\mathscr{H}^r_k}$ . }
$it_{\mathscr{H}^r_k} \leftarrow it_{\mathscr{H}^r_k} +1$\\
$d_{\mathscr{H}^r_k} \leftarrow \frac{d_{\mathscr{H}^r_k} (it_{\mathscr{H}^r_k}-1) + d(s,s^r)}{it_{\mathscr{H}^r_k}}$\\
\If{$it_{\mathscr{H}^r_k}=\gamma$}
{
    $\omega_{\mathscr{H}^r_k} \leftarrow \frac{\omega_{\mathscr{H}^r_k}d_{\beta}}{d_{\mathscr{H}^r_k}}$\\
    $\omega_{\mathscr{H}^r_k} \leftarrow \min\{size,\max\{1,\omega_{\mathscr{H}^r_k}\}\}$\\
    $it_{\mathscr{H}^r_k}, d_{\mathscr{H}^r_k} \leftarrow 0$\\
}
\caption{Updating $\omega_{\mathscr{H}^r_k}$}
\label{updateOmega}
\end{algorithm}

The input data of Algorithm~\ref{updateOmega} is both current local optimum and reference solution;  the parameter degree parameter, $\omega_{\mathscr{H}^r_k}$, to be updated; the number of iterations, $it_{\mathscr{H}^r_k}$, that the perturbation heuristic $\mathscr{H}^r_k$ was applied to a solution in AILS; the average distance between the local optima of  solutions found by perturbation heuristic $\mathscr{H}^r_k$ and their reference solutions,  $d_{\mathscr{H}^r_k}$; a fixed value $\gamma$ that represents the number of iterations that $\omega_{\mathscr{H}^r_k}$ remains with the same value;  the reference distance $d_{\beta}$, which is the ``ideal'' expected distance from a local optimum and its reference solution; and $size$, which is the number of elements of a solution. In our case, i.e., the CVRP, a solution is $\mathscr{R}$ and $|\cup_{i=1}^m R^s_i|  = |V|=n$, which means that $size$ must be $n$.

 This mechanism allows a uniform diversity control to each one of the perturbation heuristics. That is, all heuristics will have the same degree of perturbation, since the reference distance $d_{\beta}$ will be the same.

\subsubsection{Acceptation Criterion}
\label{criterioAceitacao}

The acceptance criterion establishes the rules for the reference solution $s^r$ to be updated by the current solution $s$. In this paper, we calculate a threshold called $\bar{b}$ that sets the minimum quality that the current solution must have in order to update the reference solution. This threshold is calculated according to the average quality of the solutions obtained after the local search, which is called $\bar{f}$. In addition, the threshold considers another estimator, called $\underline{f}$, which is the best solution found over the last $\min\{it,\gamma\}$ iterations,  where $it$ is the current iteration of the algorithm. To calculate $\bar{f}$, we use the weighted average of the solutions where each new value has a weight  $1/\gamma$, if the algorithm has run a number of iterations higher than or equal to $\gamma$, as described in Equation \eqref{fbar}.

\begin{align}
\bar{f} = \left \{ \begin{array}{ll}
\ \bar{f} (1- \frac{1}{\gamma}) + \frac{f(s)}{\gamma}, & \textrm{if  $ it > \gamma$}  \\
 \frac{ \bar{f} (it-1) + f(s)}{it} , & \textrm{if $ it \leq \gamma$} 
\end{array}\right. & \label{fbar}
\end{align}

The computation of $\bar{b}$ depends on a parameter $\eta\in [0,1]$, and can be calculated by the following equation: $\bar{b} = \underline{f}+\eta (\bar{f}-\underline{f})$. Therefore, $\bar{b} \in [\bar{f},\underline{f}]$. Therefore, the  control of the acceptance criterion occurs by varying  $\eta$. The higher the value of $\eta$, the greater the number of solutions accepted by the method.

In this criterion, the threshold $ \bar{b} $ is adjusted so that the average flow of accepted solutions is defined by the user. The flow is represented by  parameter $ \kappa \in [0,1] $ and indicates the percentage of solutions that are accepted. In this criterion, the value of $ \eta $ is dynamically adjusted by induction for each $ \gamma $ solutions accepted. The adjustment equation is represented by $ \eta = \max \{\epsilon,  \kappa \eta / \kappa^r \} $, where $ \kappa^r $ represents the actual percentage of solutions accepted since the last update. Constant $ \epsilon $ represents a very small value and is used to guarantee that $ \eta> 0 $.

Appendix A shows the AILS to the CVRP at length. Next section presents the hybridization of AILS with Path-Relinking proposed in this paper o solve the CVRP.

\section{AILS-PR to  the CVRP}
\label{section:ILSPR}

To enable a thorough intensification of the search space to AILS to outperform state-of-the-art CVRP metaheuristics, it was hybridized with Path-Relinking (PR). PR is an intensification search strategy  proposed by \cite{Glover1996} in combination with a tabu search heuristic \citep{Glover1998,Glover1999}. The main idea behind PR is the investigation of paths between two solutions in order to find better quality intermediate solutions.

For this, a set of elite solutions found by the proposed method, named AILS-PR, is maintained and used for the construction of PR paths. In this sense, in a given iteration, a path between the current local optimum solution and an elite solution is constructed. Then, according to a distance metric, a copy of one of the solutions of the pair, called the initial solution, is modified to generate a solution closer to the other solution of the pair, called guide or final solution.

Next sections present the algorithms for constructing the elite set and the PR procedure. 

\subsection{Elite set}
\label{updateElite}

For the PR phase of AILS-PR, family of sets $\mathscr{E}$ with elite solutions  is kept during the search process. Each set $\mathscr{E}_m$ of the family has elite solutions with the same number of routes, that is, $m$.  Moreover, the set cardinalities must not exceed the maximum limit imposed by the parameter $\sigma$. Besides, all solutions belonging to a given set must have the same number of routes. The strategy adopted to update the elite set ensures that all solutions from the sets have a minimum mutual distance of $d_{\beta}$. Algorithm \ref{Alg:updateElite} presents how the  process of updating a set $\mathscr{E}_m$ works.

\begin{algorithm}[!htb]
\LinesNumbered
\SetAlgoLined
\KwData{Solution $s$ with $m$ routes, $\mathscr{E}$}
\KwResult{Updated set $\mathscr{E}$}
Consider $MD^s=\{s_e \in \mathscr{E}_m ~|~ f(s_e) \geq f(s) ~and~ d(s,s_e)\leq d_{\beta} \}$\\
Consider  $\dot{s_e}=\arg \max_{s_e \in \mathscr{E}_m} f(s_e)$\\
Consider  $s_e^f=\arg \min^+_{s_e \in \mathscr{E}_m} (f(s_e)-f(s))$\\
\If{$|\mathscr{E}_m|<\sigma$ and ($d(s,s_e)>d_{\beta} ~\forall~s_e \in \mathscr{E}_m$)}
{
    $\mathscr{E}_m \leftarrow \mathscr{E}_m \cup \{s\}$\\
}
\ElseIf{$f(s) \leq f(\dot{s_e})$ and ($d(s,s_e)>d_{\beta} ~\forall~s_e \in \mathscr{E}_m ~|~ f(s_e) < f(s) $)}
{
    \If{$|MD^s|>0$}
    {
        $\mathscr{E}_m \leftarrow (\mathscr{E}_m \setminus MD^s) \cup \{s\}$\\
    }
    \Else
    {
        $\mathscr{E}_m \leftarrow (\mathscr{E}_m \setminus \{s_e^f\}) \cup \{s\}$\\
    }
}
\caption{Updating of elite set $\mathscr{E}$.}
\label{Alg:updateElite}
\end{algorithm}

In the case that the size of $\mathscr{E}_m$ is smaller than $\sigma$, the requirement is that the solution $s$ has a minimum distance of $d_{\beta}$ in relation to the solutions that already belong to $\mathscr{E}_m$. When $\mathscr{E}_m$ is at its cardinality limit, we check whether the quality of the candidate solution $s$ is better than or equal to the worst solution of  $\mathscr{E}_m$, called $\dot{s_e}$.  If this condition is satisfied, we verify if there is any solution in $\mathscr{E}_m$ with a better objective function value than $s$ and distance lower than $d_{\beta}$. If this solution does not exist, then $s$ may join $\mathscr{E}_m$. In order for the minimum distance criterion to be met, it is necessary to remove from the $\mathscr{E}_m$ all solutions whose distances are lower  than $d_{\beta}$ with respect to $s$. These solutions are included in the set $MD^s$. If $MD^s$ is empty, then  $s$ will replace the solution whose objective function value is the closest to $s$:  $s_e^f=\arg \min^+_{s_e \in \mathscr{E}_m} f(s_e)-f(s)$. 

\subsection{Path-Relinking}
\label{PR}


Path-Relinking (PR)  is a strategy that consists of constructing paths between pairs of solutions -- starting at an \textit{initial solution} and finishing at a \textit{guide solution} -- to find better quality intermediate solutions \citep{Glover1996}.

In the proposed AILS-PR, between a solution $s$ found after a local search and a solution of the elite family of sets $\mathscr{E}$ is constructed a path. The PR algorithm, introduced in this paper, only considers paths between pairs of solutions with the same number of routes. That is because it defines a  bijector function  $\phi: \{1,\ldots,m\} \rightarrow \{1,\ldots,m\}$ to match routes from the initial solution $s_i$ with routes of the guide solution $s_g$. A route $ R^{s_i}_k $  is matched with a route $ R^{s_g}_l$ if, and only if, $\phi(k)=l$. For this reason, we decomposed  $\mathscr{E}$ into sets of solutions with the same number of routes, that is, $\mathscr{E}=\{\mathscr{E}_{\underline{m}},\mathscr{E}_{\underline{m}+1}, \ldots , \mathscr{E}_{n}\}$. A peculiarity of AILS-PR is the similarity measure adopted to compare a pair of solutions. First, the routes of the two solutions are matched. Then, the sum of the number of vertices that appear in the matched routes describes the similarity between the solutions. Thus,  the final solution built by PR may not be the same as the guide solution, since the position of the vertices in the routes are not regarded in the similarity criterion. Algorithm~\ref{alg:PR} describes the PR  introduced in this paper.
\begin{algorithm}[!htb]
\LinesNumbered
\SetAlgoLined
\KwData{Solution $s$ with $m$ routes}
\KwResult{The best solution found $s^b$}
Randomly choose a solution $s^e \in \mathscr{E}_m$\\
Randomly choose which solution will be the initial solution  $s_i$ and which will be the guide solution $s_g$ between $s^e$ and $s$\\
Find the matching $\phi$ between the routes from $s_i$ and $s_g$ using Algorithm \ref{Alg:emparelhamento}\\ 
Make $NF=\cup_{k=1}^m R^{s_i}_k \setminus (R^{s_i}_k  \cap R^{s_g}_{\phi(k)})$\\ 
$s^b \leftarrow s_i$\\
Choose randomly a combination $C \in \mathscr{C}$\\
\While{$|NF|>0$}
{
     Calculate the priority $p_v ~\forall~ v \in NF$ according to the criterion $C$\\ 
     Let $\hat{v} \in NF$ be the vertex with the highest priority and suppose that $\hat{v} \in R^{s_g}_l$\\
     Move $\hat{v}$ for $R^{s_i}_{\phi^{-1}(l)}$ at the position that results in the minimum cost solution to $s_i$ according to Equation~\eqref{bestPosition}\\
     $NF \leftarrow NF \setminus \{\hat{v}\}$\\
    \If{$s_i$ is feasible and $f(s_i) < f(s^b)$}
    {
        $s^b \leftarrow s_i$\\
    }
}
Apply the local search to $s^b$\\
Update $\mathscr{E}$ considering solution $s^b$\\
\caption{Path-Relinking}
\label{alg:PR}
\end{algorithm}

The algorithm receives a solution $ s $ with $ m $ routes as input, and returns the best feasible solution $ s ^ b $  found in the PR search. PR will be constructed between solution $ s $  and a solution $s^e \in \mathscr{E}_m$ randomly chosen. After defining the two solutions to build the path, it is necessary to define which of them will be the initial solution, $ s_i $, and the guide solution, $ s_g $. Then, the algorithm matches  the routes in $ s_i $  with the routes in $ s_g $ using the heuristic described in Section \ref{emparelhamento}, through function $\phi$. 

The solutions in  PR are defined by changing the vertices in $ R^{s_i}_k $ to other routes that are not in  $ R^{s_g}_{\phi(k)}$, so that the similarity between $ s_i $ and $ s_g $ increases. To perform this process, we define $ NF = \cup_{k = 1}^mR^{s_i}_k - R^ {s_g}_{\phi(k)}) $ as the set with all vertices of $ s_i $ that are in $ R^{s_i}_k $ but not in $ R^{s_g}_{\phi(k)} $. The length of the path to be taken is represented by $ | NF | $. After building  $ NF $, the algorithm set the best solution of PR $ s^b $ found so far as $ s_i $. Then, we randomly choose a criterion $ C \in \mathscr{C} $ that will be used to calculate the priority of the vertices $ v \in NF $. This calculation is described in detail in Appendix B. After choosing the criterion $ C $, the path between $ s_i $ and $ s_g $ is defined by deciding the order of the vertices from $s_i$ to move to other routes in order to get closer to $s_g$.

To choose the order of the vertices in solution $ s_i $ which will be assigned to other routes in order to reach $s_g$, it is necessary to calculate the priority of each vertex in $NF$. The highest priority vertex is chosen, called $ \hat{v} $. In case of a tie, we choose the vertex that has the lowest cost regarding movement $ \mathbf{SHIFT} $, among those with the highest priority. Suppose $ \hat{v} \in R^{s_g}_l $. Therefore,  vertex $ \hat{v} $ will be moved to the route $ R^{s_i}_{\phi^{- 1} (l)} $ in the position $\hat{i}$ that results in the lowest cost for $ s_i $ according to  Equation \eqref{bestPosition}. 
\begin{align}
\hat{i}=\arg\min_{i\in\{0,1,\ldots,m_{\hat{j}}\}} d(v_i^{\hat{j}},v)+d(v_{i+1}^{\hat{j}},v)-d(v_i^{\hat{j}},v_{i+1}^{\hat{j}}) & \label{bestPosition}
\end{align}

After moving  $ \hat {v} $, we update  $ NF $, and check if solution $ s_i $  is feasible and if its cost is less than the best solution objective value. If so, we assign to $ s^b $ the solution $ s_i $. After moving all vertices in $NF$, the local search, explained in Appendix A.3, is applied to $s^b$ and the elite set is updated.

\subsection{Pairing of Routes}
\label{emparelhamento}
The pairing of the routes is necessary for the PR, as it is the reference base for assessing the similarity between pairs of solutions. This procedure consists of matching the routes according to the similarity of the vertex sets. We evaluated the use of the Hungarian algorithm that returns the optimal matching, however, the evaluated implementation presented a high computational cost. For this reason, we adopted a simple heuristic to perform this task, considering that this pairing is not preponderant for PR.  Algorithm \ref{Alg:emparelhamento} describes the route pairing procedure. 

\begin{algorithm}[!htb]
\LinesNumbered
\SetAlgoLined
\KwData{Solution $s_i$ and $s_g$.}
\KwResult{The bijector function $\phi$ defined for all $m$ indexes of the routes.}
Let $\mathscr{R}^{s_i}_e$ and $\mathscr{R}^{s_g}_e$ be the solution's route sets $s_i$ and $s_g$, respectively.\\
\For{$i = 1, \ldots, m $ }
{
    Find the route pair $(R^{s_i}_{\hat{k}},R^{s_g}_{\hat{l}})$ such that $(R^{s_i}_{\hat{k}},R^{s_g}_{\hat{l}}) = \arg \max_{ R^{s_i}_{k} \in \mathscr{R}^{s_i}_e e R^{s_g}_{l} \in \mathscr{R}^{s_g}_e} |R^{s_i}_k \cap R^{s_g}_{l}|$.\\
    $\mathscr{R}^{s_i}_e \leftarrow \mathscr{R}^{s_i}_e \setminus R^{s_i}_{\hat{k}}$ \\
    $\mathscr{R}^{s_g}_e \leftarrow \mathscr{R}^{s_i}_e \setminus R^{s_g}_{\hat{l}}$ \\
    $\phi(\hat{k}) \leftarrow \hat{l}$\\
}
\caption{Pairing $s_i$ solution routes with solution routes $s_g$.}
\label{Alg:emparelhamento}
\end{algorithm}

The method adopted in this paper iteratively matches routes from the guide solution with the initial solution. In a given iteration, it chooses the one route from each solution that has not been paired yet and that has the largest number of common vertices. This procedure repeats until all routes have been paired.

Algorithm \ref{ILSPR} describes the steps of AILS-PR to solve the CVRP.

\begin{algorithm}[!htb]
\LinesNumbered
\SetAlgoLined
\KwData{Instance data}
\KwResult{The best solution found $s^*$}
$s \leftarrow$ Construct initial solution\\ 
$s^r,s^*  \leftarrow$ Local Search($s$)\\ 
$it\leftarrow 1$\\
\Repeat{$it\geq it_{max}$}
{ 
    $s \leftarrow$ Perturbation Procedure($s^r$)\\
	 $s \leftarrow$ Local Search($s$) \\
	 $\mathcal{E}\leftarrow$ Updating of elite set  $\mathcal{E}$($s$,$\mathcal{E}$)\\
	Update the diversity control parameter $\omega_{\mathscr{H}^r_k}$ considering the distance between $s$ and $s^r$\\
    $s^r\leftarrow $ Apply acceptation criterion to $s$\\
	Update the acceptance criterion\\
	Assign $s$ to $s^*$ if $ f(s) < f(s^*)$\\
    $s^b \leftarrow$ Apply  the Path-Relinking between $s$ and $\mathcal{E}$ according to Algorithm~\ref{alg:PR}\\
    \If{$ f(s^b) < f(s^*)$}
    {
        $s^*, s^r \leftarrow s^b$
    }
    $it \leftarrow it+1$\\
    }
\caption{AILS-PR}
\label{ILSPR}
\end{algorithm}

\section{Computational Experiments}
\label{section:EC}

The computational experiments were carried out on an  Intel Xeon E5-2680v2 processor with 2.8 GHz, 10 cores and 128 GB DDR3 1866 MHz RAM. The implementation of AILS-PR was developed in Java. In the experiments, we investigate the performance of AILS-PR using the set of benchmark instances presented in \citep{Uchoa2017}. The benchmark set contains 100 CVRP instances whose sizes range from 100 to 1000 vertices. 

\subsection{Parameter tuning}
 In summary, AILS-PR has five parameters:

\begin{itemize}
    \item $\gamma$: which determines the amount of data required for the algorithm to perform an update. It is also used to calculate the probability of varying the number of routes.
    
    \item $\kappa$: which  indicates the reference flow of solutions that must be accepted to replace the reference solution.
     
    \item $d_{\beta}$: which indicates the reference distance between the reference solution and the current solution obtained after the local search. It also indicates the minimum distance between solutions from the elite set.
    
    \item $\sigma$: which sets the maximum size that a set $\mathscr{E}_m \in \mathscr{E}$ may have.
    
    \item $\varphi$: which sets the maximum size that a set $\delta(v_k)$ may have.

\end{itemize}

The fine-tuning of the AILS-PR parameters was performed using the irace package \citep{irace2016} using the 100 CVRP benchmark instances introduced by \cite{Uchoa2017}. The stop criterion considered in the tuning was 50,000 iterations without improvement. Table~\ref{parametros} shows, for each  parameter, its type,  interval of values tested by irace, and value returned by the method as the best parameter value. In the next experiments we employed such values to set the parameters.

\begin{table}[!ht]
\caption{Values employed in AILS-PR after the set up of parameters defined by irace.}
\label{parametros}
\centering
\begin{tabular}{ccccc}
\hline
\hline
Parameter & Type & Interval & Value \\
\rowcolor{Gray}
$\gamma$ & Integer & [10, 70] & 20\\
$\kappa$ & Real & [0.1, 0.6] & 0.35\\
\rowcolor{Gray}
$d_{\beta}$ & Integer & [10, 50] & 24\\
$\sigma$ & Integer & [10, 100] & 63\\
\rowcolor{Gray}
$\varphi$ & Integer & [10, 100] & 60\\
\hline
\hline
\end{tabular}
\end{table}

\subsection{Assessment of the AILS-PR components}

In this experiment, we aim at assessing the impact of  PR and the Local Search (LS)  on the solution quality of the solutions obtained by AILS-PR. For this, we carried out experiments combining the insertion and removal of these phases of the method in the algorithm. Table \ref{vercoes} presents the four tested versions of the method with and without LS and PR. On the one hand, the versions without local search -- Versions 1 and 3 -- had lines 16   from Algorithm \ref{alg:PR} and lines 2 and 6 from Algorithm~\ref{ILSPR} removed. On the other, Versions 1 and 2 -- absent of the PR -- did not considered  lines 7, 12, 13, 14 and  in Algorithm \ref {ILSPR}.

\begin{table}[!ht]
\caption{Versions of AILS-PR.}
\label{vercoes}
\centering
\begin{tabular}{lccccc}
\hline
\hline
 & & \multicolumn{4}{c}{Versions} \\ \cline{3-6}
Method    &&1  &2  &3  &4  \\
\hline
\hline
\rowcolor{Gray}
LS  &&   & \checkmark  &   & \checkmark  \\
PR  &&   &   & \checkmark  & \checkmark  \\
\hline
\hline
\end{tabular}
\end{table}

In this experiment, we ran each  version of the method 50 times to solve all 100 instances proposed in \citep{Uchoa2017}. The stop criterion was 50,000 iterations without improvement. We report the  $gap$ values between the average solutions of the 50 runs obtained by the versions  (Avg) and the best known solution (BKS) reported in the literature. The $gap$ value is calculated according to Equation~\eqref{gap}.

\begin{equation}\label{gap}
      gap=100(f-\mbox{BKS})/\mbox{BKS}
\end{equation}

\noindent where $f$ is the  average solution (Avg) in this analysis.

Similar to the classification of the instances introduced in \cite{Christiaens2020}, we divided the instances into three groups for a better analysis of the performance regarding the number of nodes. Set X-n101-247 contains instances with the number of nodes between 100 to 250 (a total of 32 instances); X-n251-491 has instances with 251 to 500 nodes (a total of 36 instances); and set X-n502-1001 contains instances with 501 to 1000 nodes (a total of 32 instances).


\begin{figure}[h!t]
\center

\definecolor{RYB1}{RGB}{230,97,1}
\definecolor{RYB2}{RGB}{200,150,250}
\pgfplotscreateplotcyclelist{colorbrewer-RYB}{
{blue!60!white,fill=blue!30!white,line width=1.5pt},
{red!70!white,fill=red!30!white,line width=1.5pt},
{orange!70!white,fill=orange!30!white,line width=1.5pt},
{violet!70!white,fill=violet!30!white,line width=1.5pt}}
\begin{tikzpicture}
\pgfplotsset{every axis legend/.append style={
at={(0.8,0.05)},
anchor=south}
},
\begin{axis}[
width=12cm,
height=7.5cm,
boxplot/draw direction=y,
cycle list name=colorbrewer-RYB,
enlarge y limits,
ymajorgrids,
ylabel={$gap$},
xtick={2.5, 7.5, 12.5},
xtick style={draw=none},
xticklabels={X-n101-247, X-n251-491, X-n502-1001},
 legend style={at={(0.5,1.15)}, anchor=north,legend columns=-1},
    legend image code/.code={
    \draw[#1, draw=none] (0cm,-0.1cm) rectangle (0.6cm,0.1cm);
    },
]
\addplot+[
boxplot prepared={
draw position=1,
lower whisker=0.0015, lower quartile=0.0800,
median=0.2381, 
upper quartile=0.4033, upper whisker=1.6625,
},
]
coordinates {};
\addlegendentry{Version 1\;};
\addplot+[
boxplot prepared={
draw position=2,
lower whisker=0.0000, lower quartile=0.0073,
median=0.0561,
upper quartile=0.0983, upper whisker=0.4696,
},
]
coordinates {};
\addlegendentry{Version 2\;};
\addplot+[
boxplot prepared={
draw position=3,
lower whisker=0.0000, lower quartile=0.0000,
median=0.0015,
upper quartile=0.0055, upper whisker=0.2302,
},
]
coordinates {};
\addlegendentry{Version 3\;};
\addplot+[
boxplot prepared={
draw position=4,
lower whisker=0.0000, lower quartile=0.0000,
median=0.0000,
upper quartile=0.0002, upper whisker=0.1017,
},
]
coordinates {};
\addlegendentry{Version 4\;};
\addplot+[
boxplot prepared={
draw position=6,
lower whisker=0.2303, lower quartile=0.5028,
median=0.7479,
upper quartile=1.0010, upper whisker=2.1823,
},
]
coordinates {};
\addplot+[
boxplot prepared={
draw position=7,
lower whisker=0.0256, lower quartile=0.2268,
median=0.2675,
upper quartile=0.3882, upper whisker=0.5093,
},
]
coordinates {};
\addplot+[
boxplot prepared={
draw position=8,
lower whisker=0.0004, lower quartile=0.0736,
median=0.1322,
upper quartile=0.2045, upper whisker=0.4990,
},
]
coordinates {};
\addplot+[
boxplot prepared={
draw position=9,
lower whisker=0.0000, lower quartile=0.0235,
median=0.0543,
upper quartile=0.0746, upper whisker=0.1891,
},
]
coordinates {};
\addplot+[
boxplot prepared={
draw position=11,
lower whisker=0.3025, lower quartile=0.6123,
median=0.8637,
upper quartile=1.3831, upper whisker=2.7679,
},
]
coordinates {};
\addplot+[
boxplot prepared={
draw position=12,
lower whisker=0.0153, lower quartile=0.243,
median=0.3842,
upper quartile=0.4782, upper whisker=0.7458,
},
]
coordinates {};
\addplot+[
boxplot prepared={
draw position=13,
lower whisker=0.0342, lower quartile=0.1968,
median=0.2586,
upper quartile=0.4262, upper whisker=0.8608,
},
]
coordinates {};
\addplot+[
boxplot prepared={
draw position=14,
lower whisker=0.003, lower quartile=0.0815,
median=0.1094,
upper quartile=0.2046, upper whisker=0.4278,
},
]
coordinates {};
\end{axis}
\end{tikzpicture}
\caption{Box plots of the 4 Versions of  AILS-PR}
\label{fig:boxplot}
\end{figure}
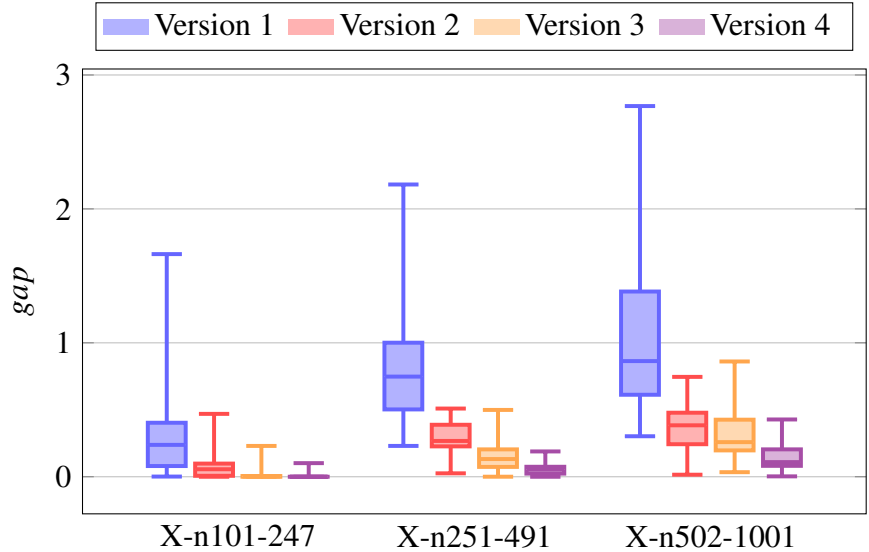

The box-plots in Figure~\ref{fig:boxplot} illustrate the results achieved by the four versions of the method. It is noteworthy that Version 4 -- AILS-PR -- was the most robust version among the methods, since its box-plot has a more condensed box, lower median, minimum and maximum values. Nevertheless, Version 3, which is AILS without the  local search, also performed pretty well, even though the upper whisker was much longer than the one observed in Version 4. On the other hand, Versions 1 and 2 had a worse performance, with longer whiskers and lengthy boxes.


\subsection{Experiment}

The set of 100 instances proposed by \cite{Uchoa2017} were tested with three algorithms recently proposed in the literature, known as UHGS \citep{Vidal2014}, ILS-SP \citep{Subramanian2013} and SISRs \citep{Christiaens2020}. 
We followed the same experimental methodology conducted by \cite{Uchoa2017}. All four solution methods were run 50 times for each instance. AILS-PR  had as stop criterion  200,000 iterations without improvement. 

Tables \ref{resultI} to \ref{resultIV} show the results of the experiment. The values reported in the tables are:

\begin{itemize}
    \item BKS: the objective function value of the best-known solutions. The values highlighted with $*$ are optimal solutions according to the information available in the repository \url{http://vrp.atd-lab.inf.puc-rio.br}\footnote{Most recent data as of 14 november 2020 - 18.00GMT.}.
    
    \item Avg: reports the average objective function value of the solutions obtained in  50 rounds.
    
    \item \textit{gap}: represents the  $gap$ between the values of the Avg column and the BKS, and is calculated according to Equation \eqref{gap}.

    \item Best: objective function value of the best solution obtained in  50 runs.
    
    \item Time: represents the average time in minutes that the algorithm took to find the best solution in each run.
\end{itemize}

The background of the top results for each instance is highlighted in dark gray; the second best results in medium gray; and the third best results in light gray. The configuration of the computer where the experiments of the reference algorithms were run are provided in the end of Table~\ref{resultIV}.

\begin{landscape}
\begin{table}[!ht]
\centering
\scalefont{0.6}
\caption{Results of experiments with instances proposed in \cite{Uchoa2017}.}
\label{resultI}
\centering
\begin{tabular}{ccccccccccccccccc}
\hline
&& \multicolumn{3}{c}{ILS-SP \citep{Subramanian2013}} &&  \multicolumn{3}{c}{UHGS \citep{Vidal2014}} && \multicolumn{3}{c}{SISRs \citep{Christiaens2020}} && \multicolumn{3}{c}{AILS-PR}\\ 
\cline{3-5} \cline{7-9} \cline{11-13} \cline{15-17}
Instance&BKS&Avg~(\textit{gap})&Best&Time&&Avg&Best&Time&&Avg~(\textit{gap})&Best&Time &&Avg~(\textit{gap})&Best&Time\\
\hline
X-n101-k25&$\mathbf{27591^*}$&\cellcolor{Gray1}27591.0~(0.0000)&\cellcolor{Gray1}27591&\cellcolor{Gray2}0.13&&\cellcolor{Gray1}27591.0~(0.0000)&\cellcolor{Gray1}27591&\cellcolor{Gray4}1.43&&\cellcolor{Gray1}27591.0~(0.0000)&\cellcolor{Gray1}27591&\cellcolor{Gray3}0.8&&\cellcolor{Gray1}27591.0~(0.0000)&\cellcolor{Gray1}27591&\cellcolor{Gray1}0.08\\
X-n106-k14&$\mathbf{26362^*}$&\cellcolor{Gray2}26375.9~(0.0527)&\cellcolor{Gray1}26362&\cellcolor{Gray3}2.01&&\cellcolor{Gray4}26381.8~(0.0751)&\cellcolor{Gray4}26378&\cellcolor{Gray4}4.04&&\cellcolor{Gray3}26381.5~(0.0740)&\cellcolor{Gray1}26362&\cellcolor{Gray2}1.3&&\cellcolor{Gray1}26362.0~(0.0000)&\cellcolor{Gray1}26362&\cellcolor{Gray1}0.84\\
X-n110-k13&$\mathbf{14971^*}$&\cellcolor{Gray1}14971.0~(0.0000)&\cellcolor{Gray1}14971&\cellcolor{Gray2}0.20&&\cellcolor{Gray1}14971.0~(0.0000)&\cellcolor{Gray1}14971&\cellcolor{Gray4}1.58&&\cellcolor{Gray4}14971.1~(0.0007)&\cellcolor{Gray1}14971&\cellcolor{Gray3}1.0&&\cellcolor{Gray1}14971.0~(0.0000)&\cellcolor{Gray1}14971&\cellcolor{Gray1}0.02\\
X-n115-k10&$\mathbf{12747^*}$&\cellcolor{Gray1}12747.0~(0.0000)&\cellcolor{Gray1}12747&\cellcolor{Gray2}0.18&&\cellcolor{Gray1}12747.0~(0.0000)&\cellcolor{Gray1}12747&\cellcolor{Gray4}1.81&&\cellcolor{Gray1}12747.0~(0.0000)&\cellcolor{Gray1}12747&\cellcolor{Gray3}0.2&&\cellcolor{Gray1}12747.0~(0.0000)&\cellcolor{Gray1}12747&\cellcolor{Gray1}0.07\\
X-n120-k6&$\mathbf{13332^*}$&\cellcolor{Gray4}13337.6~(0.0420)&\cellcolor{Gray1}13332&\cellcolor{Gray3}1.69&&\cellcolor{Gray1}13332.0~(0.0000)&\cellcolor{Gray1}13332&\cellcolor{Gray4}2.31&&\cellcolor{Gray1}13332.0~(0.0000)&\cellcolor{Gray1}13332&\cellcolor{Gray2}1.6&&\cellcolor{Gray1}13332.0~(0.0000)&\cellcolor{Gray1}13332&\cellcolor{Gray1}0.11\\
X-n125-k30&$\mathbf{55539^*}$&\cellcolor{Gray4}55673.8~(0.2427)&\cellcolor{Gray1}55539&\cellcolor{Gray2}1.43&&\cellcolor{Gray2}55542.1~(0.0056)&\cellcolor{Gray1}55539&\cellcolor{Gray3}2.66&&\cellcolor{Gray3}55556.3~(0.0311)&\cellcolor{Gray4}55542&\cellcolor{Gray4}3.1&&\cellcolor{Gray1}55539.0~(0.0000)&\cellcolor{Gray1}55539&\cellcolor{Gray1}0.73\\
X-n129-k18&$\mathbf{28940^*}$&\cellcolor{Gray4}28998.0~(0.2004)&\cellcolor{Gray4}28948&\cellcolor{Gray2}1.92&&\cellcolor{Gray2}28948.5~(0.0294)&\cellcolor{Gray1}28940&\cellcolor{Gray4}2.71&&\cellcolor{Gray3}28948.8~(0.0304)&\cellcolor{Gray1}28940&\cellcolor{Gray1}1.5&&\cellcolor{Gray1}28940.3~(0.0011)&\cellcolor{Gray1}28940&\cellcolor{Gray3}2.30\\
X-n134-k13&$\mathbf{10916^*}$&\cellcolor{Gray4}10947.4~(0.2877)&\cellcolor{Gray1}10916&\cellcolor{Gray2}2.07&&\cellcolor{Gray2}10934.9~(0.1731)&\cellcolor{Gray1}10916&\cellcolor{Gray4}3.32&&\cellcolor{Gray3}10940.1~(0.2208)&\cellcolor{Gray1}10916&\cellcolor{Gray3}2.8&&\cellcolor{Gray1}10916.0~(0.0000)&\cellcolor{Gray1}10916&\cellcolor{Gray1}0.22\\
X-n139-k10&$\mathbf{13590^*}$&\cellcolor{Gray4}13603.1~(0.0964)&\cellcolor{Gray1}13590&\cellcolor{Gray2}1.60&&\cellcolor{Gray1}13590.0~(0.0000)&\cellcolor{Gray1}13590&\cellcolor{Gray4}2.28&&\cellcolor{Gray3}13595.4~(0.0397)&\cellcolor{Gray1}13590&\cellcolor{Gray3}2.0&&\cellcolor{Gray1}13590.0~(0.0000)&\cellcolor{Gray1}13590&\cellcolor{Gray1}0.30\\
X-n143-k7&$\mathbf{15700^*}$&\cellcolor{Gray4}15745.2~(0.2879)&\cellcolor{Gray4}15726&\cellcolor{Gray1}1.64&&\cellcolor{Gray2}15700.2~(0.0013)&\cellcolor{Gray1}15700&\cellcolor{Gray4}3.10&&\cellcolor{Gray3}15705.8~(0.0369)&\cellcolor{Gray1}15700&\cellcolor{Gray2}2.1&&\cellcolor{Gray1}15700.0~(0.0000)&\cellcolor{Gray1}15700&\cellcolor{Gray3}2.57\\
X-n148-k46&$\mathbf{43448^*}$&\cellcolor{Gray3}43452.1~(0.0094)&\cellcolor{Gray1}43448&\cellcolor{Gray2}0.84&&\cellcolor{Gray1}43448.0~(0.0000)&\cellcolor{Gray1}43448&\cellcolor{Gray4}3.18&&\cellcolor{Gray4}43469.2~(0.0488)&\cellcolor{Gray1}43448&\cellcolor{Gray3}2.8&&\cellcolor{Gray1}43448.0~(0.0000)&\cellcolor{Gray1}43448&\cellcolor{Gray1}0.16\\
X-n153-k22&$\mathbf{21220^*}$&\cellcolor{Gray4}21400.0~(0.8483)&\cellcolor{Gray4}21340&\cellcolor{Gray1}0.49&&\cellcolor{Gray2}21226.3~(0.0297)&\cellcolor{Gray1}21220&\cellcolor{Gray3}5.47&&\cellcolor{Gray3}21229.4~(0.0443)&\cellcolor{Gray1}21220&\cellcolor{Gray4}5.6&&\cellcolor{Gray1}21220.0~(0.0000)&\cellcolor{Gray1}21220&\cellcolor{Gray2}4.07\\
X-n157-k13&$\mathbf{16876^*}$&\cellcolor{Gray1}16876.0~(0.0000)&\cellcolor{Gray1}16876&\cellcolor{Gray2}0.76&&\cellcolor{Gray1}16876.0~(0.0000)&\cellcolor{Gray1}16876&\cellcolor{Gray3}3.19&&\cellcolor{Gray4}16878.6~(0.0154)&\cellcolor{Gray1}16876&\cellcolor{Gray4}3.7&&\cellcolor{Gray1}16876.0~(0.0000)&\cellcolor{Gray1}16876&\cellcolor{Gray1}0.13\\
X-n162-k11&$\mathbf{14138^*}$&\cellcolor{Gray4}14160.1~(0.1563)&\cellcolor{Gray1}14138&\cellcolor{Gray2}0.54&&\cellcolor{Gray2}14141.3~(0.0233)&\cellcolor{Gray1}14138&\cellcolor{Gray3}3.32&&\cellcolor{Gray3}14157.1~(0.1351)&\cellcolor{Gray1}14138&\cellcolor{Gray4}3.4&&\cellcolor{Gray1}14138.0~(0.0000)&\cellcolor{Gray1}14138&\cellcolor{Gray1}0.18\\
X-n167-k10&$\mathbf{20557^*}$&\cellcolor{Gray4}20608.7~(0.2515)&\cellcolor{Gray4}20562&\cellcolor{Gray2}0.86&&\cellcolor{Gray3}20563.2~(0.0302)&\cellcolor{Gray1}20557&\cellcolor{Gray4}3.73&&\cellcolor{Gray2}20560.8~(0.0185)&\cellcolor{Gray1}20557&\cellcolor{Gray3}3.2&&\cellcolor{Gray1}20557.0~(0.0000)&\cellcolor{Gray1}20557&\cellcolor{Gray1}0.56\\
X-n172-k51&$\mathbf{45607^*}$&\cellcolor{Gray3}45616.1~(0.0200)&\cellcolor{Gray1}45607&\cellcolor{Gray1}0.64&&\cellcolor{Gray1}45607.0~(0.0000)&\cellcolor{Gray1}45607&\cellcolor{Gray3}3.83&&\cellcolor{Gray4}45619.2~(0.0268)&\cellcolor{Gray1}45607&\cellcolor{Gray4}5.3&&\cellcolor{Gray1}45607.0~(0.0000)&\cellcolor{Gray1}45607&\cellcolor{Gray2}0.64\\
X-n176-k26&$\mathbf{47812^*}$&\cellcolor{Gray4}48249.8~(0.9157)&\cellcolor{Gray4}48140&\cellcolor{Gray2}1.11&&\cellcolor{Gray3}47957.2~(0.3037)&\cellcolor{Gray1}47812&\cellcolor{Gray4}7.56&&\cellcolor{Gray2}47849.6~(0.0786)&\cellcolor{Gray1}47812&\cellcolor{Gray3}5.2&&\cellcolor{Gray1}47812.0~(0.0000)&\cellcolor{Gray1}47812&\cellcolor{Gray1}0.45\\
X-n181-k23&$\mathbf{25569^*}$&\cellcolor{Gray2}25571.5~(0.0098)&\cellcolor{Gray1}25569&\cellcolor{Gray1}1.59&&\cellcolor{Gray4}25591.1~(0.0864)&\cellcolor{Gray1}25569&\cellcolor{Gray4}6.28&&\cellcolor{Gray3}25579.8~(0.0422)&\cellcolor{Gray1}25569&\cellcolor{Gray3}5.5&&\cellcolor{Gray1}25569.0~(0.0002)&\cellcolor{Gray1}25569&\cellcolor{Gray2}4.26\\
X-n186-k15&$\mathbf{24145^*}$&\cellcolor{Gray4}24186.0~(0.1698)&\cellcolor{Gray1}24145&\cellcolor{Gray1}1.72&&\cellcolor{Gray2}24147.2~(0.0091)&\cellcolor{Gray1}24145&\cellcolor{Gray4}5.92&&\cellcolor{Gray3}24178.4~(0.1383)&\cellcolor{Gray4}24149&\cellcolor{Gray2}4.0&&\cellcolor{Gray1}24145.1~(0.0003)&\cellcolor{Gray1}24145&\cellcolor{Gray3}5.32\\
X-n190-k8&$\mathbf{16980^*}$&\cellcolor{Gray4}17143.1~(0.9605)&\cellcolor{Gray4}17085&\cellcolor{Gray1}2.10&&\cellcolor{Gray3}16987.9~(0.0465)&\cellcolor{Gray1}16980&\cellcolor{Gray4}12.08&&\cellcolor{Gray2}16984.9~(0.0289)&\cellcolor{Gray1}16980&\cellcolor{Gray2}9.1&&\cellcolor{Gray1}16980.3~(0.0020)&\cellcolor{Gray1}16980&\cellcolor{Gray3}9.12\\
X-n195-k51&$\mathbf{44225^*}$&\cellcolor{Gray2}44234.3~(0.0210)&\cellcolor{Gray1}44225&\cellcolor{Gray1}0.87&&\cellcolor{Gray3}44244.1~(0.0432)&\cellcolor{Gray1}44225&\cellcolor{Gray3}6.10&&\cellcolor{Gray4}44298.5~(0.1662)&\cellcolor{Gray4}44241&\cellcolor{Gray3}6.1&&\cellcolor{Gray1}44225.0~(0.0000)&\cellcolor{Gray1}44225&\cellcolor{Gray2}1.55\\
X-n200-k36&$\mathbf{58578^*}$&\cellcolor{Gray4}58697.2~(0.2035)&\cellcolor{Gray4}58626&\cellcolor{Gray3}7.48&&\cellcolor{Gray2}58626.4~(0.0826)&\cellcolor{Gray1}58578&\cellcolor{Gray4}7.97&&\cellcolor{Gray3}58636.1~(0.0992)&\cellcolor{Gray1}58578&\cellcolor{Gray2}6.7&&\cellcolor{Gray1}58578.0~(0.0000)&\cellcolor{Gray1}58578&\cellcolor{Gray1}3.05\\
X-n204-k19&$\mathbf{19565^*}$&\cellcolor{Gray3}19625.2~(0.3077)&\cellcolor{Gray4}19570&\cellcolor{Gray1}1.08&&\cellcolor{Gray2}19571.5~(0.0332)&\cellcolor{Gray1}19565&\cellcolor{Gray4}5.35&&\cellcolor{Gray4}19662.3~(0.4973)&\cellcolor{Gray1}19565&\cellcolor{Gray3}4.9&&\cellcolor{Gray1}19565.0~(0.0000)&\cellcolor{Gray1}19565&\cellcolor{Gray2}2.38\\
X-n209-k16&$\mathbf{30656^*}$&\cellcolor{Gray4}30765.4~(0.3569)&\cellcolor{Gray4}30667&\cellcolor{Gray2}3.80&&\cellcolor{Gray3}30680.4~(0.0796)&\cellcolor{Gray1}30656&\cellcolor{Gray4}8.62&&\cellcolor{Gray2}30669.4~(0.0437)&\cellcolor{Gray1}30656&\cellcolor{Gray3}6.0&&\cellcolor{Gray1}30656.0~(0.0000)&\cellcolor{Gray1}30656&\cellcolor{Gray1}1.62\\
X-n214-k11&$\mathbf{10856^*}$&\cellcolor{Gray4}11126.9~(2.4954)&\cellcolor{Gray4}10985&\cellcolor{Gray1}2.26&&\cellcolor{Gray2}10877.4~(0.1971)&\cellcolor{Gray1}10856&\cellcolor{Gray4}10.22&&\cellcolor{Gray3}10908.6~(0.4845)&\cellcolor{Gray3}10873&\cellcolor{Gray2}8.7&&\cellcolor{Gray1}10867.0~(0.1017)&\cellcolor{Gray2}10857&\cellcolor{Gray3}9.80\\
X-n219-k73&$\mathbf{117595^*}$&\cellcolor{Gray1}117595.0~(0.0000)&\cellcolor{Gray1}117595&\cellcolor{Gray2}0.85&&\cellcolor{Gray3}117604.9~(0.0084)&\cellcolor{Gray1}117595&\cellcolor{Gray3}7.73&&\cellcolor{Gray4}117650.4~(0.0471)&\cellcolor{Gray1}117595&\cellcolor{Gray4}8.0&&\cellcolor{Gray1}117595.0~(0.0000)&\cellcolor{Gray1}117595&\cellcolor{Gray1}0.31\\
\hline
\end{tabular}
\end{table}
\end{landscape}

\begin{landscape}
\begin{table}[!ht]
\centering
\scalefont{0.6}
\caption{Results of experiments with instances proposed in \cite{Uchoa2017}.}
\label{resultII}
\centering
\begin{tabular}{crccccccccccccccc}
\hline
  && \multicolumn{3}{c}{ILS-SP \citep{Subramanian2013}} &&  \multicolumn{3}{c}{UHGS \citep{Vidal2014}} && \multicolumn{3}{c}{SISRs \citep{Christiaens2020}} && \multicolumn{3}{c}{AILS-PR}\\ 
\cline{3-5} \cline{7-9} \cline{11-13} \cline{15-17}
Instance&BKS&Avg~(\textit{gap})&Best&Time&&Avg~(\textit{gap})&Best&Time&&Avg~(\textit{gap})&Best&Time &&Avg~(\textit{gap})&Best&Time\\
\hline
X-n223-k34&$\mathbf{40437^*}$&\cellcolor{Gray4}40533.5~(0.2386)&\cellcolor{Gray4}40471&\cellcolor{Gray4}8.48&&\cellcolor{Gray2}40499.0~(0.1533)&\cellcolor{Gray1}40437&\cellcolor{Gray3}8.26&&\cellcolor{Gray3}40529.9~(0.2297)&\cellcolor{Gray3}40448&\cellcolor{Gray2}7.6&&\cellcolor{Gray1}40437.2~(0.0004)&\cellcolor{Gray1}40437&\cellcolor{Gray1}7.08\\
X-n228-k23&$\mathbf{25742^*}$&\cellcolor{Gray4}25795.8~(0.2090)&\cellcolor{Gray3}25743&\cellcolor{Gray1}2.40&&\cellcolor{Gray2}25779.3~(0.1449)&\cellcolor{Gray1}25742&\cellcolor{Gray3}9.80&&\cellcolor{Gray3}25790.9~(0.1900)&\cellcolor{Gray4}25744&\cellcolor{Gray4}10.5&&\cellcolor{Gray1}25742.9~(0.0037)&\cellcolor{Gray1}25742&\cellcolor{Gray2}2.66\\
X-n233-k16&$\mathbf{19230^*}$&\cellcolor{Gray4}19336.7~(0.5549)&\cellcolor{Gray4}19266&\cellcolor{Gray1}3.01&&\cellcolor{Gray3}19288.4~(0.3037)&\cellcolor{Gray1}19230&\cellcolor{Gray2}6.84&&\cellcolor{Gray2}19269.7~(0.2064)&\cellcolor{Gray3}19232&\cellcolor{Gray4}8.1&&\cellcolor{Gray1}19230.1~(0.0004)&\cellcolor{Gray1}19230&\cellcolor{Gray3}7.10\\
X-n237-k14&$\mathbf{27042^*}$&\cellcolor{Gray3}27078.8~(0.1361)&\cellcolor{Gray1}27042&\cellcolor{Gray2}3.46&&\cellcolor{Gray2}27067.3~(0.0936)&\cellcolor{Gray1}27042&\cellcolor{Gray4}8.90&&\cellcolor{Gray4}27089.7~(0.1764)&\cellcolor{Gray1}27042&\cellcolor{Gray3}7.2&&\cellcolor{Gray1}27042.0~(0.0000)&\cellcolor{Gray1}27042&\cellcolor{Gray1}1.98\\
X-n242-k48&$\mathbf{82751^*}$&\cellcolor{Gray2}82874.2~(0.1489)&\cellcolor{Gray2}82774&\cellcolor{Gray4}17.83&&\cellcolor{Gray4}82948.7~(0.2389)&\cellcolor{Gray4}82804&\cellcolor{Gray2}12.42&&\cellcolor{Gray3}82884.4~(0.1612)&\cellcolor{Gray3}82775&\cellcolor{Gray1}9.9&&\cellcolor{Gray1}82808.1~(0.0690)&\cellcolor{Gray1}82764&\cellcolor{Gray3}13.12\\
X-n247-k50&$\mathbf{37274^*}$&\cellcolor{Gray4}37507.2~(0.6256)&\cellcolor{Gray4}37289&\cellcolor{Gray2}2.06&&\cellcolor{Gray2}37284.4~(0.0279)&\cellcolor{Gray1}37274&\cellcolor{Gray4}20.41&&\cellcolor{Gray3}37323.2~(0.1320)&\cellcolor{Gray1}37274&\cellcolor{Gray3}18.4&&\cellcolor{Gray1}37274.0~(0.0000)&\cellcolor{Gray1}37274&\cellcolor{Gray1}1.17\\
X-n251-k28&$\mathbf{38684^*}$&\cellcolor{Gray4}38840.0~(0.4033)&\cellcolor{Gray4}38727&\cellcolor{Gray3}10.77&&\cellcolor{Gray3}38796.4~(0.2906)&\cellcolor{Gray3}38699&\cellcolor{Gray4}11.69&&\cellcolor{Gray2}38791.0~(0.2766)&\cellcolor{Gray2}38687&\cellcolor{Gray2}9.8&&\cellcolor{Gray1}38713.7~(0.0768)&\cellcolor{Gray1}38684&\cellcolor{Gray1}9.50\\
X-n256-k16&$18839$&\cellcolor{Gray3}18883.9~(0.2383)&\cellcolor{Gray2}18880&\cellcolor{Gray2}2.02&&\cellcolor{Gray2}18880.0~(0.2176)&\cellcolor{Gray2}18880&\cellcolor{Gray3}6.52&&\cellcolor{Gray4}18888.9~(0.2649)&\cellcolor{Gray2}18880&\cellcolor{Gray4}11.5&&\cellcolor{Gray1}18874.6~(0.1891)&\cellcolor{Gray1}18851&\cellcolor{Gray1}1.65\\
X-n261-k13&$26558$&\cellcolor{Gray4}26869.0~(1.1710)&\cellcolor{Gray4}26706&\cellcolor{Gray1}6.67&&\cellcolor{Gray2}26629.6~(0.2696)&\cellcolor{Gray1}26558&\cellcolor{Gray4}12.67&&\cellcolor{Gray3}26642.3~(0.3174)&\cellcolor{Gray1}26558&\cellcolor{Gray3}11.8&&\cellcolor{Gray1}26558.4~(0.0014)&\cellcolor{Gray1}26558&\cellcolor{Gray2}10.75\\
X-n266-k58&$\mathbf{75478^*}$&\cellcolor{Gray2}75563.3~(0.1130)&\cellcolor{Gray1}75478&\cellcolor{Gray1}10.03&&\cellcolor{Gray4}75759.3~(0.3727)&\cellcolor{Gray4}75517&\cellcolor{Gray4}21.36&&\cellcolor{Gray3}75617.8~(0.1852)&\cellcolor{Gray1}75478&\cellcolor{Gray2}10.8&&\cellcolor{Gray1}75533.5~(0.0735)&\cellcolor{Gray1}75478&\cellcolor{Gray3}12.79\\
X-n270-k35&$\mathbf{35291^*}$&\cellcolor{Gray3}35363.4~(0.2052)&\cellcolor{Gray4}35324&\cellcolor{Gray1}9.07&&\cellcolor{Gray4}35367.2~(0.2159)&\cellcolor{Gray2}35303&\cellcolor{Gray3}11.25&&\cellcolor{Gray2}35362.2~(0.2018)&\cellcolor{Gray3}35323&\cellcolor{Gray4}11.4&&\cellcolor{Gray1}35304.2~(0.0374)&\cellcolor{Gray1}35291&\cellcolor{Gray2}10.30\\
X-n275-k28&$\mathbf{21245^*}$&\cellcolor{Gray2}21256.0~(0.0518)&\cellcolor{Gray1}21245&\cellcolor{Gray2}3.59&&\cellcolor{Gray4}21280.6~(0.1676)&\cellcolor{Gray1}21245&\cellcolor{Gray3}12.04&&\cellcolor{Gray3}21268.6~(0.1111)&\cellcolor{Gray1}21245&\cellcolor{Gray4}13.3&&\cellcolor{Gray1}21245.0~(0.0000)&\cellcolor{Gray1}21245&\cellcolor{Gray1}3.11\\
X-n280-k17&$33503$&\cellcolor{Gray4}33769.4~(0.7952)&\cellcolor{Gray4}33624&\cellcolor{Gray1}9.62&&\cellcolor{Gray2}33605.8~(0.3068)&\cellcolor{Gray2}33505&\cellcolor{Gray4}19.09&&\cellcolor{Gray3}33628.1~(0.3734)&\cellcolor{Gray3}33529&\cellcolor{Gray3}17.7&&\cellcolor{Gray1}33526.4~(0.0699)&\cellcolor{Gray1}33503&\cellcolor{Gray2}13.51\\
X-n284-k15&$\mathbf{20215^*}$&\cellcolor{Gray4}20448.5~(1.1551)&\cellcolor{Gray4}20295&\cellcolor{Gray1}8.64&&\cellcolor{Gray2}20286.4~(0.3532)&\cellcolor{Gray1}20227&\cellcolor{Gray4}19.91&&\cellcolor{Gray3}20286.6~(0.3542)&\cellcolor{Gray3}20240&\cellcolor{Gray2}15.3&&\cellcolor{Gray1}20247.1~(0.1586)&\cellcolor{Gray2}20229&\cellcolor{Gray3}18.20\\
X-n289-k60&$\mathbf{95151^*}$&\cellcolor{Gray3}95450.6~(0.3149)&\cellcolor{Gray4}95315&\cellcolor{Gray2}16.11&&\cellcolor{Gray4}95469.5~(0.3347)&\cellcolor{Gray3}95244&\cellcolor{Gray4}21.28&&\cellcolor{Gray2}95352.2~(0.2115)&\cellcolor{Gray2}95233&\cellcolor{Gray1}14.3&&\cellcolor{Gray1}95205.9~(0.0577)&\cellcolor{Gray1}95151&\cellcolor{Gray3}17.81\\
X-n294-k50&$47161$&\cellcolor{Gray2}47254.7~(0.1987)&\cellcolor{Gray3}47190&\cellcolor{Gray1}12.42&&\cellcolor{Gray3}47259.0~(0.2078)&\cellcolor{Gray2}47171&\cellcolor{Gray2}14.70&&\cellcolor{Gray4}47274.5~(0.2407)&\cellcolor{Gray4}47210&\cellcolor{Gray2}14.7&&\cellcolor{Gray1}47175.9~(0.0316)&\cellcolor{Gray1}47167&\cellcolor{Gray4}17.22\\
X-n298-k31&$\mathbf{34231^*}$&\cellcolor{Gray4}34356.0~(0.3652)&\cellcolor{Gray4}34239&\cellcolor{Gray2}6.92&&\cellcolor{Gray3}34292.1~(0.1785)&\cellcolor{Gray1}34231&\cellcolor{Gray3}10.93&&\cellcolor{Gray2}34276.0~(0.1315)&\cellcolor{Gray3}34234&\cellcolor{Gray4}14.5&&\cellcolor{Gray1}34231.0~(0.0000)&\cellcolor{Gray1}34231&\cellcolor{Gray1}3.93\\
X-n303-k21&$21736$&\cellcolor{Gray4}21895.8~(0.7352)&\cellcolor{Gray4}21812&\cellcolor{Gray1}14.15&&\cellcolor{Gray3}21850.9~(0.5286)&\cellcolor{Gray2}21748&\cellcolor{Gray3}17.28&&\cellcolor{Gray2}21776.5~(0.1863)&\cellcolor{Gray3}21751&\cellcolor{Gray4}17.3&&\cellcolor{Gray1}21748.6~(0.0578)&\cellcolor{Gray1}21738&\cellcolor{Gray2}15.48\\
X-n308-k13&$25859$&\cellcolor{Gray3}26101.1~(0.9362)&\cellcolor{Gray3}25901&\cellcolor{Gray1}9.53&&\cellcolor{Gray2}25895.4~(0.1408)&\cellcolor{Gray1}25859&\cellcolor{Gray2}15.31&&\cellcolor{Gray4}26207.7~(1.3485)&\cellcolor{Gray4}25931&\cellcolor{Gray4}25.7&&\cellcolor{Gray1}25866.8~(0.0300)&\cellcolor{Gray1}25859&\cellcolor{Gray3}24.01\\
X-n313-k71&$94043$&\cellcolor{Gray4}94297.3~(0.2704)&\cellcolor{Gray4}94192&\cellcolor{Gray1}17.50&&\cellcolor{Gray3}94265.2~(0.2363)&\cellcolor{Gray3}94093&\cellcolor{Gray4}22.41&&\cellcolor{Gray2}94182.4~(0.1482)&\cellcolor{Gray2}94063&\cellcolor{Gray2}18.9&&\cellcolor{Gray1}94065.1~(0.0235)&\cellcolor{Gray1}94044&\cellcolor{Gray3}22.14\\
X-n317-k53&$\mathbf{78355^*}$&\cellcolor{Gray1}78356.0~(0.0013)&\cellcolor{Gray1}78355&\cellcolor{Gray1}8.56&&\cellcolor{Gray2}78387.8~(0.0419)&\cellcolor{Gray1}78355&\cellcolor{Gray4}22.37&&\cellcolor{Gray4}78392.4~(0.0477)&\cellcolor{Gray1}78355&\cellcolor{Gray3}22.0&&\cellcolor{Gray3}78389.4~(0.0439)&\cellcolor{Gray1}78355&\cellcolor{Gray2}11.26\\
X-n322-k28&$\mathbf{29834^*}$&\cellcolor{Gray4}29991.3~(0.5273)&\cellcolor{Gray4}29877&\cellcolor{Gray2}14.68&&\cellcolor{Gray3}29956.1~(0.4093)&\cellcolor{Gray3}29870&\cellcolor{Gray3}15.16&&\cellcolor{Gray2}29927.6~(0.3137)&\cellcolor{Gray2}29849&\cellcolor{Gray4}16.9&&\cellcolor{Gray1}29849.1~(0.0507)&\cellcolor{Gray1}29834&\cellcolor{Gray1}14.02\\
X-n327-k20&$27532$&\cellcolor{Gray4}27812.4~(1.0185)&\cellcolor{Gray3}27599&\cellcolor{Gray3}19.13&&\cellcolor{Gray2}27628.2~(0.3494)&\cellcolor{Gray2}27564&\cellcolor{Gray2}18.19&&\cellcolor{Gray3}27631.4~(0.3610)&\cellcolor{Gray4}27608&\cellcolor{Gray4}21.6&&\cellcolor{Gray1}27572.8~(0.1480)&\cellcolor{Gray1}27543&\cellcolor{Gray1}12.33\\
X-n331-k15&$\mathbf{31102^*}$&\cellcolor{Gray4}31235.5~(0.4292)&\cellcolor{Gray3}31105&\cellcolor{Gray2}15.70&&\cellcolor{Gray3}31159.6~(0.1852)&\cellcolor{Gray2}31103&\cellcolor{Gray4}24.43&&\cellcolor{Gray2}31128.2~(0.0842)&\cellcolor{Gray4}31122&\cellcolor{Gray3}20.4&&\cellcolor{Gray1}31103.5~(0.0048)&\cellcolor{Gray1}31102&\cellcolor{Gray1}8.23\\
X-n336-k84&$139111$&\cellcolor{Gray3}139461.0~(0.2516)&\cellcolor{Gray2}139197&\cellcolor{Gray1}21.41&&\cellcolor{Gray4}139534.9~(0.3047)&\cellcolor{Gray4}139210&\cellcolor{Gray4}37.96&&\cellcolor{Gray2}139373.4~(0.1886)&\cellcolor{Gray3}139209&\cellcolor{Gray2}22.8&&\cellcolor{Gray1}139197.3~(0.0620)&\cellcolor{Gray1}139160&\cellcolor{Gray3}29.14\\
\hline
\end{tabular}
\end{table}
\end{landscape}

\begin{landscape}
\begin{table}[!ht]
\centering
\scalefont{0.6}
\caption{Results of experiments with instances proposed in \cite{Uchoa2017}.}
\label{resultIII}
\centering
\begin{tabular}{crccccccccccccccc}
\hline
  && \multicolumn{3}{c}{ILS-SP \citep{Subramanian2013}} &&  \multicolumn{3}{c}{UHGS \citep{Vidal2014}} && \multicolumn{3}{c}{SISRs \citep{Christiaens2020}} && \multicolumn{3}{c}{AILS-PR}\\ 
\cline{3-5} \cline{7-9} \cline{11-13} \cline{15-17}
Instance&BKS&Avg~(\textit{gap})&Best&Time&&Avg~(\textit{gap})&Best&Time&&Avg~(\textit{gap})&Best&Time &&Avg~(\textit{gap})&Best&Time\\
\hline
X-n344-k43&$42050$&\cellcolor{Gray4}42284.0~(0.5565)&\cellcolor{Gray4}42146&\cellcolor{Gray4}22.58&&\cellcolor{Gray3}42208.8~(0.3776)&\cellcolor{Gray3}42099&\cellcolor{Gray3}21.67&&\cellcolor{Gray2}42158.5~(0.2580)&\cellcolor{Gray2}42079&\cellcolor{Gray2}21.5&&\cellcolor{Gray1}42072.8~(0.0543)&\cellcolor{Gray1}42056&\cellcolor{Gray1}19.81\\
X-n351-k40&$25896$&\cellcolor{Gray4}26150.3~(0.9820)&\cellcolor{Gray4}26021&\cellcolor{Gray1}25.21&&\cellcolor{Gray3}26014.0~(0.4557)&\cellcolor{Gray3}25946&\cellcolor{Gray4}33.73&&\cellcolor{Gray2}25982.1~(0.3325)&\cellcolor{Gray2}25938&\cellcolor{Gray2}26.5&&\cellcolor{Gray1}25937.2~(0.1591)&\cellcolor{Gray1}25920&\cellcolor{Gray3}32.17\\
X-n359-k29&$51505$&\cellcolor{Gray4}52076.5~(1.1096)&\cellcolor{Gray4}51706&\cellcolor{Gray4}48.86&&\cellcolor{Gray3}51721.7~(0.4207)&\cellcolor{Gray3}51509&\cellcolor{Gray3}34.85&&\cellcolor{Gray2}51577.8~(0.1413)&\cellcolor{Gray1}51505&\cellcolor{Gray1}23.1&&\cellcolor{Gray1}51548.0~(0.0834)&\cellcolor{Gray2}51507&\cellcolor{Gray2}25.41\\
X-n367-k17&$22814$&\cellcolor{Gray4}23003.2~(0.8293)&\cellcolor{Gray4}22902&\cellcolor{Gray2}13.13&&\cellcolor{Gray3}22838.4~(0.1070)&\cellcolor{Gray1}22814&\cellcolor{Gray3}22.02&&\cellcolor{Gray2}22833.4~(0.0850)&\cellcolor{Gray1}22814&\cellcolor{Gray4}36.1&&\cellcolor{Gray1}22814.0~(0.0000)&\cellcolor{Gray1}22814&\cellcolor{Gray1}6.29\\
X-n376-k94&$\mathbf{147713^*}$&\cellcolor{Gray1}147713.0~(0.0000)&\cellcolor{Gray1}147713&\cellcolor{Gray2}7.10&&\cellcolor{Gray3}147750.2~(0.0252)&\cellcolor{Gray3}147717&\cellcolor{Gray3}28.26&&\cellcolor{Gray4}147783.6~(0.0478)&\cellcolor{Gray4}147721&\cellcolor{Gray4}32.0&&\cellcolor{Gray1}147713.0~(0.0000)&\cellcolor{Gray1}147713&\cellcolor{Gray1}3.40\\
X-n384-k52&$65938$&\cellcolor{Gray4}66372.5~(0.6590)&\cellcolor{Gray4}66116&\cellcolor{Gray3}34.47&&\cellcolor{Gray3}66270.2~(0.5038)&\cellcolor{Gray3}66081&\cellcolor{Gray4}40.20&&\cellcolor{Gray2}66107.4~(0.2569)&\cellcolor{Gray2}65963&\cellcolor{Gray1}25.9&&\cellcolor{Gray1}66001.9~(0.0969)&\cellcolor{Gray1}65953&\cellcolor{Gray2}28.15\\
X-n393-k38&$\mathbf{38260^*}$&\cellcolor{Gray4}38457.4~(0.5159)&\cellcolor{Gray3}38298&\cellcolor{Gray2}20.82&&\cellcolor{Gray2}38374.9~(0.3003)&\cellcolor{Gray2}38269&\cellcolor{Gray3}28.65&&\cellcolor{Gray3}38394.1~(0.3505)&\cellcolor{Gray4}38331&\cellcolor{Gray4}30.4&&\cellcolor{Gray1}38266.6~(0.0172)&\cellcolor{Gray1}38260&\cellcolor{Gray1}19.90\\
X-n401-k29&$66163$&\cellcolor{Gray4}66715.1~(0.8345)&\cellcolor{Gray4}66453&\cellcolor{Gray4}60.36&&\cellcolor{Gray3}66365.4~(0.3059)&\cellcolor{Gray3}66243&\cellcolor{Gray3}49.52&&\cellcolor{Gray2}66248.5~(0.1292)&\cellcolor{Gray2}66189&\cellcolor{Gray2}38.0&&\cellcolor{Gray1}66208.0~(0.0680)&\cellcolor{Gray1}66174&\cellcolor{Gray1}29.71\\
X-n411-k19&$19712$&\cellcolor{Gray4}19954.9~(1.2322)&\cellcolor{Gray4}19792&\cellcolor{Gray1}23.76&&\cellcolor{Gray2}19743.8~(0.1613)&\cellcolor{Gray1}19718&\cellcolor{Gray3}34.71&&\cellcolor{Gray3}19768.5~(0.2866)&\cellcolor{Gray3}19731&\cellcolor{Gray4}58.4&&\cellcolor{Gray1}19736.9~(0.1265)&\cellcolor{Gray2}19723&\cellcolor{Gray2}28.26\\
X-n420-k130&$\mathbf{107798^*}$&\cellcolor{Gray2}107838.0~(0.0371)&\cellcolor{Gray1}107798&\cellcolor{Gray1}22.19&&\cellcolor{Gray4}107924.1~(0.1170)&\cellcolor{Gray1}107798&\cellcolor{Gray4}53.19&&\cellcolor{Gray3}107879.2~(0.0753)&\cellcolor{Gray4}107817&\cellcolor{Gray3}47.9&&\cellcolor{Gray1}107826.2~(0.0262)&\cellcolor{Gray1}107798&\cellcolor{Gray2}33.30\\
X-n429-k61&$65449$&\cellcolor{Gray4}65746.6~(0.4547)&\cellcolor{Gray4}65563&\cellcolor{Gray3}38.22&&\cellcolor{Gray3}65648.5~(0.3048)&\cellcolor{Gray3}65501&\cellcolor{Gray4}41.45&&\cellcolor{Gray2}65593.6~(0.2209)&\cellcolor{Gray2}65485&\cellcolor{Gray2}35.0&&\cellcolor{Gray1}65485.8~(0.0562)&\cellcolor{Gray1}65449&\cellcolor{Gray1}25.05\\
X-n439-k37&$\mathbf{36391^*}$&\cellcolor{Gray2}36441.6~(0.1390)&\cellcolor{Gray2}36395&\cellcolor{Gray3}39.63&&\cellcolor{Gray3}36451.1~(0.1652)&\cellcolor{Gray2}36395&\cellcolor{Gray2}34.55&&\cellcolor{Gray4}36473.8~(0.2275)&\cellcolor{Gray4}36426&\cellcolor{Gray4}42.1&&\cellcolor{Gray1}36408.2~(0.0473)&\cellcolor{Gray1}36394&\cellcolor{Gray1}14.90\\
X-n449-k29&$55233$&\cellcolor{Gray4}56204.9~(1.7596)&\cellcolor{Gray4}55761&\cellcolor{Gray3}59.94&&\cellcolor{Gray3}55553.1~(0.5795)&\cellcolor{Gray3}55378&\cellcolor{Gray4}64.92&&\cellcolor{Gray2}55411.2~(0.3226)&\cellcolor{Gray2}55272&\cellcolor{Gray1}38.0&&\cellcolor{Gray1}55282.5~(0.0896)&\cellcolor{Gray1}55239&\cellcolor{Gray2}38.62\\
X-n459-k26&$24139$&\cellcolor{Gray4}24462.4~(1.3397)&\cellcolor{Gray4}24209&\cellcolor{Gray4}60.59&&\cellcolor{Gray3}24272.6~(0.5535)&\cellcolor{Gray3}24181&\cellcolor{Gray2}42.80&&\cellcolor{Gray2}24242.2~(0.4275)&\cellcolor{Gray2}24175&\cellcolor{Gray3}56.5&&\cellcolor{Gray1}24155.5~(0.0683)&\cellcolor{Gray1}24140&\cellcolor{Gray1}35.37\\
X-n469-k138&$\mathbf{221824^*}$&\cellcolor{Gray2}222182.0~(0.1614)&\cellcolor{Gray2}221909&\cellcolor{Gray1}36.32&&\cellcolor{Gray4}222617.1~(0.3575)&\cellcolor{Gray4}222070&\cellcolor{Gray4}86.65&&\cellcolor{Gray3}222227.1~(0.1817)&\cellcolor{Gray3}221984&\cellcolor{Gray2}48.0&&\cellcolor{Gray1}221916.9~(0.0419)&\cellcolor{Gray1}221835&\cellcolor{Gray3}59.30\\
X-n480-k70&$89449$&\cellcolor{Gray4}89871.2~(0.4720)&\cellcolor{Gray4}89694&\cellcolor{Gray2}50.40&&\cellcolor{Gray3}89760.1~(0.3478)&\cellcolor{Gray3}89535&\cellcolor{Gray4}66.96&&\cellcolor{Gray2}89559.2~(0.1232)&\cellcolor{Gray2}89458&\cellcolor{Gray3}50.5&&\cellcolor{Gray1}89457.4~(0.0094)&\cellcolor{Gray1}89449&\cellcolor{Gray1}38.33\\
X-n491-k59&$66483$&\cellcolor{Gray4}67226.7~(1.1186)&\cellcolor{Gray4}66965&\cellcolor{Gray2}52.23&&\cellcolor{Gray3}66898.0~(0.6242)&\cellcolor{Gray3}66633&\cellcolor{Gray4}71.94&&\cellcolor{Gray2}66645.5~(0.2444)&\cellcolor{Gray2}66517&\cellcolor{Gray1}51.4&&\cellcolor{Gray1}66534.6~(0.0776)&\cellcolor{Gray1}66485&\cellcolor{Gray3}63.62\\
X-n502-k39&$69226$&\cellcolor{Gray4}69346.8~(0.1745)&\cellcolor{Gray4}69284&\cellcolor{Gray4}80.75&&\cellcolor{Gray3}69328.8~(0.1485)&\cellcolor{Gray3}69253&\cellcolor{Gray3}63.61&&\cellcolor{Gray2}69274.7~(0.0703)&\cellcolor{Gray2}69243&\cellcolor{Gray2}60.9&&\cellcolor{Gray1}69237.1~(0.0160)&\cellcolor{Gray1}69226&\cellcolor{Gray1}35.40\\
X-n513-k21&$24201$&\cellcolor{Gray4}24434.0~(0.9628)&\cellcolor{Gray4}24332&\cellcolor{Gray3}35.04&&\cellcolor{Gray3}24296.6~(0.3950)&\cellcolor{Gray1}24201&\cellcolor{Gray2}33.09&&\cellcolor{Gray2}24292.1~(0.3764)&\cellcolor{Gray3}24238&\cellcolor{Gray4}77.1&&\cellcolor{Gray1}24208.8~(0.0324)&\cellcolor{Gray1}24201&\cellcolor{Gray1}31.58\\
X-n524-k153&$\mathbf{154593^*}$&\cellcolor{Gray4}155005.0~(0.2665)&\cellcolor{Gray3}154709&\cellcolor{Gray1}27.27&&\cellcolor{Gray3}154979.5~(0.2500)&\cellcolor{Gray4}154774&\cellcolor{Gray3}80.70&&\cellcolor{Gray2}154807.2~(0.1386)&\cellcolor{Gray2}154651&\cellcolor{Gray4}151.4&&\cellcolor{Gray1}154597.7~(0.0030)&\cellcolor{Gray1}154593&\cellcolor{Gray2}29.13\\
X-n536-k96&$94846$&\cellcolor{Gray4}95700.7~(0.9011)&\cellcolor{Gray4}95524&\cellcolor{Gray1}62.07&&\cellcolor{Gray3}95330.6~(0.5109)&\cellcolor{Gray3}95122&\cellcolor{Gray4}107.53&&\cellcolor{Gray2}95173.2~(0.3450)&\cellcolor{Gray2}95006&\cellcolor{Gray2}74.7&&\cellcolor{Gray1}94934.2~(0.0930)&\cellcolor{Gray1}94889&\cellcolor{Gray3}85.92\\
X-n548-k50&$\mathbf{86700^*}$&\cellcolor{Gray3}86874.1~(0.2008)&\cellcolor{Gray2}86710&\cellcolor{Gray2}63.95&&\cellcolor{Gray4}86998.5~(0.3443)&\cellcolor{Gray4}86822&\cellcolor{Gray4}84.24&&\cellcolor{Gray2}86798.0~(0.1130)&\cellcolor{Gray2}86710&\cellcolor{Gray3}64.5&&\cellcolor{Gray1}86770.7~(0.0815)&\cellcolor{Gray1}86700&\cellcolor{Gray1}34.75\\
X-n561-k42&$42717$&\cellcolor{Gray4}43131.3~(0.9699)&\cellcolor{Gray4}42952&\cellcolor{Gray3}68.86&&\cellcolor{Gray2}42866.4~(0.3497)&\cellcolor{Gray2}42756&\cellcolor{Gray2}60.60&&\cellcolor{Gray3}42868.1~(0.3537)&\cellcolor{Gray3}42774&\cellcolor{Gray4}73.8&&\cellcolor{Gray1}42747.5~(0.0715)&\cellcolor{Gray1}42717&\cellcolor{Gray1}59.09\\
X-n573-k30&$50673$&\cellcolor{Gray4}51173.0~(0.9867)&\cellcolor{Gray4}51092&\cellcolor{Gray2}112.03&&\cellcolor{Gray3}50915.1~(0.4778)&\cellcolor{Gray3}50780&\cellcolor{Gray4}188.15&&\cellcolor{Gray2}50804.6~(0.2597)&\cellcolor{Gray2}50737&\cellcolor{Gray3}113.0&&\cellcolor{Gray1}50728.4~(0.1094)&\cellcolor{Gray1}50681&\cellcolor{Gray1}90.38\\
X-n586-k159&$190316$&\cellcolor{Gray4}190919.0~(0.3168)&\cellcolor{Gray4}190612&\cellcolor{Gray1}78.54&&\cellcolor{Gray3}190838.0~(0.2743)&\cellcolor{Gray3}190543&\cellcolor{Gray4}175.29&&\cellcolor{Gray2}190600.7~(0.1496)&\cellcolor{Gray2}190484&\cellcolor{Gray2}86.3&&\cellcolor{Gray1}190369.6~(0.0282)&\cellcolor{Gray1}190324&\cellcolor{Gray3}100.79\\
X-n599-k92&$108451$&\cellcolor{Gray4}109384.0~(0.8603)&\cellcolor{Gray4}109056&\cellcolor{Gray1}72.96&&\cellcolor{Gray3}109064.2~(0.5654)&\cellcolor{Gray3}108813&\cellcolor{Gray4}125.91&&\cellcolor{Gray2}108688.6~(0.2191)&\cellcolor{Gray2}108548&\cellcolor{Gray2}75.4&&\cellcolor{Gray1}108561.7~(0.1020)&\cellcolor{Gray1}108465&\cellcolor{Gray3}79.55\\
\hline
\end{tabular}
\end{table}
\end{landscape}

\begin{landscape}
\begin{table}[!ht]
\centering
\scalefont{0.6}
\caption{Results of experiments with instances proposed in \cite{Uchoa2017}.}
\label{resultIV}
\centering
\begin{threeparttable}
\begin{tabular}{crccccccccccccccc}
\hline
  && \multicolumn{3}{c}{ILS-SP \citep{Subramanian2013}} &&  \multicolumn{3}{c}{UHGS \citep{Vidal2014}} && \multicolumn{3}{c}{SISRs \citep{Christiaens2020}} && \multicolumn{3}{c}{AILS-PR}\\ 
\cline{3-5} \cline{7-9} \cline{11-13} \cline{15-17}
Instance&BKS&Avg~(\textit{gap})&Best&Time\tnote{1} &&Avg~(\textit{gap})&Best&Time\tnote{1}&&Avg~(\textit{gap})&Best&Time\tnote{2} &&Avg~(\textit{gap})&Best&Time\tnote{3}\\
\hline
X-n613-k62&$59535$&\cellcolor{Gray4}60444.2~(1.5272)&\cellcolor{Gray4}60229&\cellcolor{Gray1}74.80&&\cellcolor{Gray3}59960.0~(0.7139)&\cellcolor{Gray3}59778&\cellcolor{Gray4}117.31&&\cellcolor{Gray2}59731.3~(0.3297)&\cellcolor{Gray2}59585&\cellcolor{Gray3}88.1&&\cellcolor{Gray1}59600.5~(0.1100)&\cellcolor{Gray1}59536&\cellcolor{Gray2}87.19\\
X-n627-k43&$62164$&\cellcolor{Gray4}62905.6~(1.1930)&\cellcolor{Gray4}62783&\cellcolor{Gray3}162.67&&\cellcolor{Gray3}62524.1~(0.5793)&\cellcolor{Gray3}62366&\cellcolor{Gray4}239.68&&\cellcolor{Gray2}62317.1~(0.2463)&\cellcolor{Gray2}62219&\cellcolor{Gray2}89.3&&\cellcolor{Gray1}62303.2~(0.2239)&\cellcolor{Gray1}62215&\cellcolor{Gray1}71.40\\
X-n641-k35&$63692$&\cellcolor{Gray4}64606.1~(1.4352)&\cellcolor{Gray4}64462&\cellcolor{Gray3}140.42&&\cellcolor{Gray3}64192.0~(0.7850)&\cellcolor{Gray3}63839&\cellcolor{Gray4}158.81&&\cellcolor{Gray2}63850.3~(0.2485)&\cellcolor{Gray2}63750&\cellcolor{Gray2}92.5&&\cellcolor{Gray1}63786.4~(0.1482)&\cellcolor{Gray1}63692&\cellcolor{Gray1}89.73\\
X-n655-k131&$\mathbf{106780^*}$&\cellcolor{Gray1}106782.0~(0.0019)&\cellcolor{Gray1}106780&\cellcolor{Gray1}47.24&&\cellcolor{Gray4}106899.1~(0.1115)&\cellcolor{Gray4}106829&\cellcolor{Gray4}150.48&&\cellcolor{Gray3}106844.6~(0.0605)&\cellcolor{Gray3}106813&\cellcolor{Gray3}109.6&&\cellcolor{Gray2}106787.8~(0.0073)&\cellcolor{Gray1}106780&\cellcolor{Gray2}63.72\\
X-n670-k130&$146332$&\cellcolor{Gray4}147676.0~(0.9185)&\cellcolor{Gray4}147045&\cellcolor{Gray1}61.24&&\cellcolor{Gray3}147222.7~(0.6087)&\cellcolor{Gray3}146705&\cellcolor{Gray4}264.10&&\cellcolor{Gray1}146720.4~(0.2654)&\cellcolor{Gray2}146451&\cellcolor{Gray3}198.9&&\cellcolor{Gray2}146833.6~(0.3428)&\cellcolor{Gray1}146406&\cellcolor{Gray2}121.36\\
X-n685-k75&$68205$&\cellcolor{Gray4}68988.2~(1.1483)&\cellcolor{Gray4}68646&\cellcolor{Gray1}73.85&&\cellcolor{Gray3}68654.1~(0.6585)&\cellcolor{Gray3}68425&\cellcolor{Gray4}156.71&&\cellcolor{Gray2}68369.0~(0.2405)&\cellcolor{Gray2}68271&\cellcolor{Gray3}135.1&&\cellcolor{Gray1}68297.8~(0.1361)&\cellcolor{Gray1}68236&\cellcolor{Gray2}85.71\\
X-n701-k44&$81923$&\cellcolor{Gray4}83042.2~(1.3662)&\cellcolor{Gray4}82888&\cellcolor{Gray3}210.08&&\cellcolor{Gray3}82487.4~(0.6889)&\cellcolor{Gray3}82293&\cellcolor{Gray4}253.17&&\cellcolor{Gray2}82065.4~(0.1738)&\cellcolor{Gray2}81974&\cellcolor{Gray2}122.5&&\cellcolor{Gray1}82049.8~(0.1548)&\cellcolor{Gray1}81941&\cellcolor{Gray1}90.86\\
X-n716-k35&$43379$&\cellcolor{Gray4}44171.6~(1.8272)&\cellcolor{Gray4}44021&\cellcolor{Gray3}225.79&&\cellcolor{Gray3}43641.4~(0.6049)&\cellcolor{Gray3}43525&\cellcolor{Gray4}264.28&&\cellcolor{Gray2}43483.8~(0.2416)&\cellcolor{Gray2}43426&\cellcolor{Gray2}158.3&&\cellcolor{Gray1}43423.2~(0.1020)&\cellcolor{Gray1}43379&\cellcolor{Gray1}128.18\\
X-n733-k159&$136187$&\cellcolor{Gray4}137045.0~(0.6300)&\cellcolor{Gray4}136832&\cellcolor{Gray1}111.56&&\cellcolor{Gray3}136587.6~(0.2942)&\cellcolor{Gray3}136366&\cellcolor{Gray4}244.53&&\cellcolor{Gray2}136389.3~(0.1485)&\cellcolor{Gray2}136255&\cellcolor{Gray3}143.2&&\cellcolor{Gray1}136300.9~(0.0836)&\cellcolor{Gray1}136223&\cellcolor{Gray2}130.44\\
X-n749-k98&$77269$&\cellcolor{Gray4}78275.9~(1.3031)&\cellcolor{Gray4}77952&\cellcolor{Gray1}127.24&&\cellcolor{Gray3}77864.9~(0.7712)&\cellcolor{Gray3}77715&\cellcolor{Gray4}313.88&&\cellcolor{Gray2}77509.2~(0.3109)&\cellcolor{Gray2}77380&\cellcolor{Gray3}146.3&&\cellcolor{Gray1}77450.8~(0.2353)&\cellcolor{Gray1}77332&\cellcolor{Gray2}127.49\\
X-n766-k71&$114418$&\cellcolor{Gray4}115738.0~(1.1537)&\cellcolor{Gray4}115443&\cellcolor{Gray3}242.11&&\cellcolor{Gray3}115147.9~(0.6379)&\cellcolor{Gray3}114683&\cellcolor{Gray4}382.99&&\cellcolor{Gray2}114761.1~(0.2999)&\cellcolor{Gray2}114590&\cellcolor{Gray2}174.4&&\cellcolor{Gray1}114619.9~(0.1764)&\cellcolor{Gray1}114453&\cellcolor{Gray1}152.01\\
X-n783-k48&$72393$&\cellcolor{Gray4}73722.9~(1.8371)&\cellcolor{Gray4}73447&\cellcolor{Gray2}235.48&&\cellcolor{Gray3}73009.6~(0.8517)&\cellcolor{Gray3}72781&\cellcolor{Gray4}269.70&&\cellcolor{Gray1}72660.7~(0.3698)&\cellcolor{Gray2}72492&\cellcolor{Gray1}170.2&&\cellcolor{Gray2}72666.1~(0.3773)&\cellcolor{Gray1}72470&\cellcolor{Gray3}238.52\\
X-n801-k40&$73305$&\cellcolor{Gray4}74005.7~(0.9559)&\cellcolor{Gray4}73830&\cellcolor{Gray4}432.64&&\cellcolor{Gray3}73731.0~(0.5811)&\cellcolor{Gray3}73587&\cellcolor{Gray3}289.24&&\cellcolor{Gray1}73436.7~(0.1797)&\cellcolor{Gray2}73347&\cellcolor{Gray2}137.1&&\cellcolor{Gray2}73478.1~(0.2362)&\cellcolor{Gray1}73320&\cellcolor{Gray1}92.87\\
X-n819-k171&$158121$&\cellcolor{Gray4}159425.0~(0.8247)&\cellcolor{Gray4}159164&\cellcolor{Gray1}148.91&&\cellcolor{Gray3}158899.3~(0.4922)&\cellcolor{Gray3}158611&\cellcolor{Gray4}374.28&&\cellcolor{Gray2}158423.0~(0.1910)&\cellcolor{Gray2}158305&\cellcolor{Gray2}172.5&&\cellcolor{Gray1}158261.9~(0.0891)&\cellcolor{Gray1}158147&\cellcolor{Gray3}230.68\\
X-n837-k142&$193737$&\cellcolor{Gray4}195027.0~(0.6659)&\cellcolor{Gray4}194804&\cellcolor{Gray3}173.17&&\cellcolor{Gray3}194476.5~(0.3817)&\cellcolor{Gray3}194266&\cellcolor{Gray4}463.36&&\cellcolor{Gray2}193976.9~(0.1238)&\cellcolor{Gray2}193824&\cellcolor{Gray2}166.8&&\cellcolor{Gray1}193924.1~(0.0966)&\cellcolor{Gray1}193797&\cellcolor{Gray1}111.48\\
X-n856-k95&$88965$&\cellcolor{Gray4}89277.6~(0.3514)&\cellcolor{Gray3}89060&\cellcolor{Gray2}153.65&&\cellcolor{Gray3}89238.7~(0.3076)&\cellcolor{Gray4}89118&\cellcolor{Gray4}288.43&&\cellcolor{Gray2}89131.3~(0.1869)&\cellcolor{Gray2}89050&\cellcolor{Gray3}160.0&&\cellcolor{Gray1}89040.6~(0.0850)&\cellcolor{Gray1}88968&\cellcolor{Gray1}111.16\\
X-n876-k59&$99299$&\cellcolor{Gray4}100417.0~(1.1259)&\cellcolor{Gray4}100177&\cellcolor{Gray3}409.31&&\cellcolor{Gray3}99884.1~(0.5892)&\cellcolor{Gray3}99715&\cellcolor{Gray4}495.38&&\cellcolor{Gray2}99483.2~(0.1855)&\cellcolor{Gray2}99388&\cellcolor{Gray1}217.4&&\cellcolor{Gray1}99439.8~(0.1418)&\cellcolor{Gray1}99333&\cellcolor{Gray2}225.08\\
X-n895-k37&$53860$&\cellcolor{Gray4}54958.5~(2.0395)&\cellcolor{Gray4}54713&\cellcolor{Gray4}410.17&&\cellcolor{Gray3}54439.8~(1.0765)&\cellcolor{Gray3}54172&\cellcolor{Gray3}321.89&&\cellcolor{Gray2}54085.8~(0.4192)&\cellcolor{Gray2}53993&\cellcolor{Gray2}212.5&&\cellcolor{Gray1}54028.8~(0.3134)&\cellcolor{Gray1}53906&\cellcolor{Gray1}97.50\\
X-n916-k207&$329179$&\cellcolor{Gray4}330948.0~(0.5374)&\cellcolor{Gray4}330639&\cellcolor{Gray2}226.08&&\cellcolor{Gray3}330198.3~(0.3096)&\cellcolor{Gray3}329836&\cellcolor{Gray4}560.81&&\cellcolor{Gray2}329509.5~(0.1004)&\cellcolor{Gray2}329299&\cellcolor{Gray1}215.3&&\cellcolor{Gray1}329348.9~(0.0516)&\cellcolor{Gray1}329228&\cellcolor{Gray3}227.75\\
X-n936-k151&$132715$&\cellcolor{Gray4}134530.0~(1.3676)&\cellcolor{Gray4}133592&\cellcolor{Gray1}202.50&&\cellcolor{Gray3}133512.9~(0.6012)&\cellcolor{Gray3}133140&\cellcolor{Gray4}531.50&&\cellcolor{Gray2}133117.3~(0.3031)&\cellcolor{Gray2}133014&\cellcolor{Gray3}412.7&&\cellcolor{Gray1}132986.6~(0.2046)&\cellcolor{Gray1}132839&\cellcolor{Gray2}255.76\\
X-n957-k87&$85465$&\cellcolor{Gray4}85936.6~(0.5518)&\cellcolor{Gray4}85697&\cellcolor{Gray3}311.20&&\cellcolor{Gray3}85822.6~(0.4184)&\cellcolor{Gray3}85672&\cellcolor{Gray4}432.90&&\cellcolor{Gray2}85620.0~(0.1814)&\cellcolor{Gray2}85546&\cellcolor{Gray2}202.4&&\cellcolor{Gray1}85577.6~(0.1317)&\cellcolor{Gray1}85481&\cellcolor{Gray1}151.11\\
X-n979-k58&$118976$&\cellcolor{Gray4}120253.0~(1.0733)&\cellcolor{Gray4}119994&\cellcolor{Gray4}687.22&&\cellcolor{Gray3}119502.1~(0.4422)&\cellcolor{Gray3}119194&\cellcolor{Gray3}553.96&&\cellcolor{Gray1}119120.4~(0.1214)&\cellcolor{Gray2}119065&\cellcolor{Gray1}276.6&&\cellcolor{Gray2}119496.0~(0.4371)&\cellcolor{Gray1}118987&\cellcolor{Gray2}309.22\\
X-n1001-k43&$72355$&\cellcolor{Gray4}73985.4~(2.2533)&\cellcolor{Gray4}73776&\cellcolor{Gray4}792.75&&\cellcolor{Gray3}72956.0~(0.8306)&\cellcolor{Gray3}72742&\cellcolor{Gray3}549.03&&\cellcolor{Gray1}72528.1~(0.2392)&\cellcolor{Gray2}72415&\cellcolor{Gray2}284.3&&\cellcolor{Gray2}72612.9~(0.3565)&\cellcolor{Gray1}72399&\cellcolor{Gray1}123.66\\
Average&&~(0.6296)&&71.71&&~(0.2969)&&98.79&&~(0.1969)&&60.01&&~(0.0710)&&47.10\\
\hline
\end{tabular}
\begin{tablenotes}
\item[1] 3.07 GHz Xeon CPU and 16 GB of RAM, running on Oracle Linux Server 6.4.
\item[2] Xeon E5-2650 v2 CPU at 2.60 GHz.
\item[3]Intel Xeon E5-2680v2 processor with 2.8 GHz, 10 cores and 128 GB DDR3 1866 MHz RAM.
\end{tablenotes}
\end{threeparttable}
\end{table}
\end{landscape}

According to Tables~\ref{resultI} to \ref{resultIV}, considering all instances, AILS-PR achieved the best results in the vast majority of the instances. More specifically,  it obtained the top best solutions (column `Best')  in 96\%  of the instances. Considering the average solutions (column `Avg'), it achieved the top results in 93\% of the instances. On the other hand, the best solutions of UHGS, SISRs and ILS-SP were the best overall results in 40\%, 30\% and 23\% of the instances, respectively. And the average solutions were the top results in 8\% of the instances in the literature methods.  Moreover, AILS-PR improved the best known solutions of 3 instances, already available in CVRPLIB: X-n384-k52 (65938), X-n641-k35 (63692), X-n716-k35 (43379).

Regarding the computational time, AILS-PR presented the shortest time for 51\% of instances. AILS-PR obtained the shortest average time, 47.1 minutes. SISRs, ILS-SP and UHGS presented an average time of 60.01, 71.71 and 98.79 minutes, respectively.

\subsubsection{Analysis of the results}

This section shows a summary of the results in three different ways. For ease of notation, we refer to the \textit{gaps} between average solutions of 50 runs and BKSs as \textit{average gaps}, and between the best solutions of 50 runs and BKSs as \textit{best gaps}. The first analysis shows a summary of the \textit{average and best gaps} achieved by the four solution methods  in the instances divided into  groups  X-n101-247, X-n251-491 and X-n502-1001 as well as the mean times. The second is the performance profiles \citep{Dolan2002} of the \textit{average gaps} of the solutions found by the heuristic methods. The third analysis presents a ranking of the methods per group of instances regarding the \textit{average and best gaps} and the average times in minutes that the algorithms took to find the best solution in each run.

Table~\ref{resultsummary3} presents the statistics  of the average and best results of the four heuristic methods to solve the instances divided into three groups according to their size.
 For each algorithm, we present  the minimum (`Min'), mean (`Mean') and maximum (`Max') \textit{average} and \textit{best gaps}, and average times considering the instances of the given group.  Besides, to a better visualization of the values for each algorithm, we show an interval representing such minimum, mean and maximum values.

\def\ILSSP101{(0,0,a)(1,0,b)(2,0,c)}
\def\UHGS101{(1,0)(2,0)(3,0)}
\def\SISR101{(2,0)(3,0)(4,0)}
\def\AILSPR101{(3,0)(4,0)(5,0)}

\newcommand{\plotReguaBest}[4]{%
    \begin{tikzpicture}[xscale=#1, yscale=1]
        \draw[
        blue!50!white,
        line width=1pt](0,0) -- (#2,0);
        \foreach \x in {0,#4,...,#2}
        {
            \draw[line width=0.1pt,violet!50!white](\x,0.05) to (\x ,-0.05);
        }
        \foreach \x/\txt in{ 0,#3,...,#2}
        {
            \draw[line width=1pt,red!50!white](\x,0.1) to (\x ,-0.1);
            \node at (\x,0.3) {\tiny {\txt}};
            \node at (0,0.3) {};
        }
    \end{tikzpicture}%
}
\newcommand{\plotReguagambi}{%
    \begin{tikzpicture}[xscale=0.03, yscale=1]
        \draw[
        blue!50!white,
        line width=1pt](0,0) -- (80,0);
        \foreach \x in {0,5,...,80}
        {
            \draw[line width=0.1pt,violet!50!white](\x,0.05) to (\x ,-0.05);
        }
        \foreach \x/\txt in{ 0/0,20/200,40/400,60/600,80/800}
        {
            \draw[line width=1pt,red!50!white](\x,0.1) to (\x ,-0.1);
            \node at (\x,0.3) {\tiny \txt};
            \node at (0,0.3) {};
        }
    \end{tikzpicture}%
}

\newcommand{\plot}[4]{%
    \begin{tikzpicture}[xscale=#4, yscale=1]
        \draw[
        blue!50!white,
        line width=1pt](#1,0) -- (#3,0);
        \foreach \x/\txt in{#1,#2,#3}
        {
            \draw[line width=1pt,red!50!white](\x,0.1) to (\x ,-0.1);
            \node at (\x,0.3) {};
            \node at (0,0.3) {};
        }
    \end{tikzpicture}%
}

\addtolength{\tabcolsep}{-4pt}  
\begin{table}[!ht]
\centering
\scalefont{0.62}
\caption{Summary of the results obtained by  the evaluated methods.}
\label{resultsummary3}
\centering
\begin{tabular}{l|cccl|cccl|cccl}
\hline
&\multicolumn{4}{c|}{Statistics of \textit{average gaps}}&\multicolumn{4}{c|}{Statistics of \textit{best gaps}}&\multicolumn{4}{c}{Statistics of average times}  \\ 
\hline
\multirow{2}{*}{X-n101-247}&\multirow{2}{*}{Min}&\multirow{2}{*}{Mean}&\multirow{2}{*}{Max}&\multicolumn{1}{c|}{Interval}&\multirow{2}{*}{Min}&\multirow{2}{*}{Mean}&\multirow{2}{*}{Max}&\multicolumn{1}{c|}{Interval}&\multirow{2}{*}{Min}&\multirow{2}{*}{Mean}&\multirow{2}{*}{Max}&\multicolumn{1}{c}{Interval}\\
& & & &\plotReguaBest{0.8}{2.5}{0.5}{0.1}& & &&\plotReguaBest{1.5}{1.5}{0.5}{0.1}& & &&\plotReguaBest{0.1}{20}{4}{1}\\
\hline
ILS-SP&0.00&0.31&2.50&\plot{0.00}{0.31}{2.5}{0.8} &0.00&0.12&1.19&\plot{0.00}{0.12}{1.19}{1.5}&0.13&2.41&17.83&\plot{0.13}{2.41}{17.83}{0.1}\\
UHGS&0.00&{0.07}&{0.30}&\plot{0.00}{0.07}{0.3}{0.8} &0.00&0.00&0.06&\plot{0.00}{0.00}{0.06}{1.5}&1.43&6.01&20.41&\plot{1.43}{6.01}{20.41}{0.1}\\
SISR&0.00&{0.11}&{0.50}&\plot{0.00}{0.11}{0.5}{0.8} &0.00&0.01&0.16&\plot{0.00}{0.01}{0.16}{1.5}&0.20&5.20&18.40&\plot{0.20}{5.20}{18.40}{0.1}\\
AILS-PR &0.00&{0.01}&{0.10}&\plot{0.00}{0.01}{0.1}{0.8}&0.00&0.00&0.02&\plot{0.00}{0.00}{0.02}{1.5}&0.02&2.62&13.12&\plot{0.02}{2.62}{13.12}{0.1}\\
\hline
\multirow{2}{*}{X-n251-491}&\multirow{2}{*}{Min}&\multirow{2}{*}{Mean}&\multirow{2}{*}{Max}&\multicolumn{1}{c|}{Interval}&\multirow{2}{*}{Min}&\multirow{2}{*}{Mean}&\multirow{2}{*}{Max}&\multicolumn{1}{c|}{Interval}&\multirow{2}{*}{Min}&\multirow{2}{*}{Mean}&\multirow{2}{*}{Max}&\multicolumn{1}{c}{Interval}\\
&&&&\plotReguaBest{1}{2}{0.5}{0.1}& & &&\plotReguaBest{2}{1}{0.5}{0.1}& & &&\plotReguaBest{0.025}{80}{20}{5}\\
\hline
ILS-SP& 0.00&{0.59}&{1.76}&\plot{0.00}{0.59}{1.76}{1} &0.00&0.23&0.96&\plot{0.00}{0.23}{0.96}{2}&2.02&23.12&60.59 &\plot{2.02}{23.12}{60.59}{0.025} \\
UHGS&0.00&{0.28}&{0.62}&\plot{0.03}{0.28}{0.62}{1}&0.00&0.07&0.26&\plot{0.00}{0.07}{0.26}{2}&6.52&25.70&86.65&\plot{6.52}{25.70}{86.65}{0.025}   \\
SISR&0.00&{0.25}&{1.35}&\plot{0.05}{0.25}{1.35}{1} &0.00&0.07&0.28&\plot{0.00}{0.07}{0.28}{2}&9.80&23.93&58.40&\plot{9.80}{23.93}{58.40}{0.025}\\
AILS-PR&0.00&{0.06}&{0.19}&\plot{0.00}{0.06}{0.19}{1}&0.00&0.01&0.09&\plot{0.00}{0.01}{0.09}{2}&1.65&17.51&63.62 &\plot{1.65}{17.51}{63.62}{0.025} \\
\hline
\multirow{2}{*}{X-n251-491}&\multirow{2}{*}{Min}&\multirow{2}{*}{Mean}&\multirow{2}{*}{Max}&\multicolumn{1}{c|}{Interval}&\multirow{2}{*}{Min}&\multirow{2}{*}{Mean}&\multirow{2}{*}{Max}&\multicolumn{1}{c|}{Interval}&\multirow{2}{*}{Min}&\multirow{2}{*}{Mean}&\multirow{2}{*}{Max}&\multicolumn{1}{c}{Interval}\\
&&&&\plotReguaBest{0.8}{2.5}{0.5}{0.1}& & &&\plotReguaBest{1}{2}{0.5}{0.1}& & &&\plotReguagambi\\
\hline
ILS-SP&0.00&{0.99}&{2.25} &\plot{0.00}{0.99}{2.25}{0.8}&0.00&0.72&1.96&\plot{0.00}{0.72}{1.96}{1}&27.27&195.67&792.75&\plot{2.727}{19.567}{79.275}{0.03}\\
UHGS&{0.11}&{0.52}&{1.08}&\plot{0.11}{0.52}{1.08}{0.8}&0.00&0.28&0.58&\plot{0.00}{0.28}{0.58}{1} &33.09&268.61&560.81&\plot{3.309}{26.861}{56.081}{0.03}\\
SISR&{0.06}&{0.22}&{0.42} &\plot{0.06}{0.22}{0.42}{0.8}&0.01&0.10&0.25&\plot{0.01}{0.10}{0.25}{1}&60.90&151.97&412.70&\plot{6.090}{15.197}{41.270}{0.03}\\
AILS-PR&0.00&{0.15}&{0.44}&\plot{0.00}{0.15}{0.44}{0.8}&0.00&0.03&0.11&\plot{0.00}{0.03}{0.11}{1}&29.13&120.92&309.22 &\plot{2.913}{12.092}{30.922}{0.03}\\
\hline
\end{tabular}
\end{table}
\addtolength{\tabcolsep}{4pt} 

The results in Table~\ref{resultsummary3} demonstrate the outstanding performance of AILS-PR in comparison to the literature heuristics, considering both \textit{gaps} and times. The mean \textit{average} and \textit{best gaps} are very close to the minimum \textit{gaps}, which are zero in all tested cases.

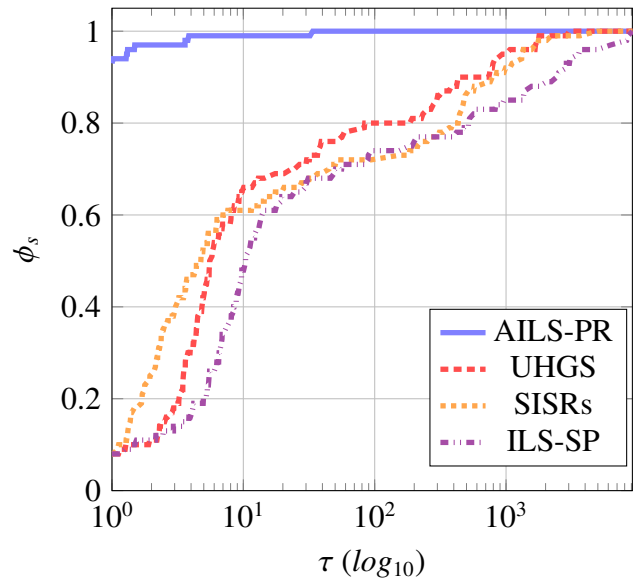
\begin{figure}[htbp]
\center
\begin{tikzpicture}
\pgfplotsset{every axis legend/.append style={
at={(0.8,0.05)},
anchor=south}
},
\begin{axis}[
width=8.5cm,
height=8cm,
xmode=log,
log basis x=10,
legend columns=1,
xmin=1,
xmax=9157,
ymin=0,
ymax=1.05,
grid=major,
xlabel={$\tau$ ($log_{10}$)},
ylabel={$\phi_s$ },
]

\addplot[/tikz/solid, blue!50!white,line width=2pt]
table[x=x,y=y] {dados/PD/ILSPRGap.txt};

\addplot[/tikz/densely dashed, red!70!white,line width=2pt]
table[x=x,y=y] {dados/PD/UHGSGap.txt};

\addplot[/tikz/dotted, orange!70!white,line width=2pt]
table[x=x,y=y] {dados/PD/ASBRRGap.txt};

\addplot[/tikz/dashdotdotted,	violet!70!white,line width=2pt]
table[x=x,y=y] {dados/PD/ILSGap.txt};

\legend{AILS-PR, UHGS, SISRs, ILS-SP }
\end{axis}
\end{tikzpicture}
\caption{Performance profile \citep{Dolan2002} comparing the \textit{average gaps} obtained by the ILS-SP, UHGS, SISRs and AILS-PR.}
\label{fig:dolanGap}
\end{figure}

In  Figure \ref{fig:dolanGap}, the \textit{average gaps} of the four algorithms were analyzed through the performance profiles. The performance profiles display the percentage of problems solved by a particular algorithm based on a performance factor $\tau \in \mathbb{R}$ that relates the analyzed algorithm to the best overall result. Thus, we have in the $y$-axis  $\phi_h(\tau)= P(r_{ph} \leq \tau ~|~1<h< n_h)$ which represents the cumulative probability that the performance ratio $r_{ph}$ associated with algorithm $h$ is within a performance factor $\tau$ indicated in the $x$-axis. This means that when $\tau$ is 1,  $\phi_h$ indicates the percentage of problems for which algorithm $h$ achieved the best results. The value of $\tau$ for $\phi_h=1$ represents the performance factor for the algorithm to achieve the best results. The values to be analyzed by the performance profiles have to be positive. For this reason, we replaced all gaps less than or equal to zero with 0.0001. In this figure, we can verify that the curve of the performance profile of AILS-PR dominates the curves of SISRs, ILS-SP and UHGS. One curve dominates the others indicates that this algorithm had a greater percentage of problems solved for any performance factor in comparison to the other algorithms. AILS-PR produced the best results in 93\% of the instances. When  $\tau$ is 33.77, which is a relatively low-performance factor, AILS-PR achieved the best results in 100\% of the instances. The second best results were obtained by UHGS. This algorithm performed better in 8\% of instances. Moreover, UHGS reaches the best result for 100\% of instances only with a performance factor of 3037. 
SISRs achieved the best results in  8\% of the instances   and, in the worst case, it presented an \textit{average gap} 4973 times greater than the best result. ILS-SP had the best result in 8\% of the instances and only with a performance factor of 9157 that it achieves the best results. 


Table \ref{resultsummary2} shows the third analysis of the results achieved by the algorithms, by ranking them according to their performance in solving  the instances proposed in \cite{Uchoa2017}. This table reports the proportion of instances per group that a given algorithm was ranked in positions 1, 2, 3 and 4. Therefore, rank 1 refers  to the proportion of instances from the indicated group that the algorithm presented the best result among the 4 algorithms; rank 2 represents the proportion of instances that the algorithm achieved the second best result; and so on. To a better visualization, we highlight the background of the results with the top values in dark gray, second top in medium gray and third top in light gray.

\begin{table}[!ht]
\centering
\scalefont{0.6}
\caption{Proportion of results obtained by the heuristic methods ranked into four positions.}
\label{resultsummary2}
\centering
\begin{tabular}{ccccccccccccccccc}
\hline
& \multicolumn{3}{c}{ILS-SP} && \multicolumn{3}{c}{UHGS} && \multicolumn{3}{c}{SISRs} && \multicolumn{3}{c}{AILS-PR}\\ 
\cline{2-4} \cline{6-8} \cline{10-12}  \cline{14-16}
Rank&Avg&Best&Time&&Avg&Best&Time&&Avg&Best&Time&&Avg&Best&Time \\
\hline
X-n101-247&&&&&&&&&&&&&&&\\
Rank 1&\cellcolor{Gray3}0.16&\cellcolor{Gray4}0.53&\cellcolor{Gray2}0.34&&\cellcolor{Gray2}0.25&\cellcolor{Gray2}0.94&\cellcolor{Gray4}0.00&&\cellcolor{Gray4}0.09&\cellcolor{Gray3}0.75&\cellcolor{Gray3}0.06&&\cellcolor{Gray1}1.00&\cellcolor{Gray1}0.97&\cellcolor{Gray1}0.59\\
Rank 2&\cellcolor{Gray3}0.13&\cellcolor{Gray2}0.03&\cellcolor{Gray1}0.50&&\cellcolor{Gray1}0.44&\cellcolor{Gray3}0.00&\cellcolor{Gray4}0.06&&\cellcolor{Gray2}0.16&\cellcolor{Gray3}0.00&\cellcolor{Gray2}0.25&&\cellcolor{Gray4}0.00&\cellcolor{Gray2}0.03&\cellcolor{Gray3}0.19\\
Rank 3&\cellcolor{Gray3}0.13&\cellcolor{Gray2}0.03&\cellcolor{Gray4}0.09&&\cellcolor{Gray2}0.22&\cellcolor{Gray3}0.00&\cellcolor{Gray2}0.28&&\cellcolor{Gray1}0.50&\cellcolor{Gray1}0.13&\cellcolor{Gray1}0.44&&\cellcolor{Gray4}0.00&\cellcolor{Gray3}0.00&\cellcolor{Gray3}0.22\\
Rank 4&\cellcolor{Gray1}0.59&\cellcolor{Gray1}0.41&\cellcolor{Gray3}0.06&&\cellcolor{Gray3}0.09&\cellcolor{Gray3}0.06&\cellcolor{Gray1}0.66&&\cellcolor{Gray2}0.25&\cellcolor{Gray2}0.13&\cellcolor{Gray2}0.25&&\cellcolor{Gray4}0.00&\cellcolor{Gray4}0.00&\cellcolor{Gray4}0.00\\
X-n251-491&&&&&&&&&&&&&&&\\
Rank 1&\cellcolor{Gray2}0.06&\cellcolor{Gray4}0.14&\cellcolor{Gray2}0.42&&\cellcolor{Gray3}0.00&\cellcolor{Gray2}0.25&\cellcolor{Gray4}0.00&&\cellcolor{Gray3}0.00&\cellcolor{Gray3}0.17&\cellcolor{Gray3}0.14&&\cellcolor{Gray1}0.97&\cellcolor{Gray1}0.92&\cellcolor{Gray1}0.44\\
Rank 2&\cellcolor{Gray3}0.17&\cellcolor{Gray3}0.11&\cellcolor{Gray2}0.31&&\cellcolor{Gray2}0.25&\cellcolor{Gray2}0.25&\cellcolor{Gray4}0.14&&\cellcolor{Gray1}0.56&\cellcolor{Gray1}0.39&\cellcolor{Gray2}0.31&&\cellcolor{Gray4}0.00&\cellcolor{Gray4}0.08&\cellcolor{Gray3}0.28\\
Rank 3&\cellcolor{Gray3}0.14&\cellcolor{Gray3}0.14&\cellcolor{Gray4}0.17&&\cellcolor{Gray1}0.56&\cellcolor{Gray1}0.42&\cellcolor{Gray1}0.36&&\cellcolor{Gray2}0.28&\cellcolor{Gray2}0.22&\cellcolor{Gray3}0.19&&\cellcolor{Gray4}0.03&\cellcolor{Gray4}0.00&\cellcolor{Gray2}0.25\\
Rank 4&\cellcolor{Gray1}0.64&\cellcolor{Gray1}0.61&\cellcolor{Gray3}0.11&&\cellcolor{Gray2}0.19&\cellcolor{Gray3}0.08&\cellcolor{Gray1}0.50&&\cellcolor{Gray3}0.17&\cellcolor{Gray2}0.22&\cellcolor{Gray2}0.36&&\cellcolor{Gray4}0.00&\cellcolor{Gray4}0.00&\cellcolor{Gray4}0.03\\
X-n502-1001&&&&&&&&&&&&&&&\\
Rank 1&\cellcolor{Gray3}0.03&\cellcolor{Gray2}0.03&\cellcolor{Gray2}0.38&&\cellcolor{Gray4}0.00&\cellcolor{Gray2}0.03&\cellcolor{Gray4}0.00&&\cellcolor{Gray2}0.16&\cellcolor{Gray3}0.00&\cellcolor{Gray3}0.13&&\cellcolor{Gray1}0.81&\cellcolor{Gray1}1.00&\cellcolor{Gray1}0.50\\
Rank 2&\cellcolor{Gray4}0.00&\cellcolor{Gray2}0.03&\cellcolor{Gray3}0.16&&\cellcolor{Gray3}0.03&\cellcolor{Gray2}0.03&\cellcolor{Gray4}0.06&&\cellcolor{Gray1}0.78&\cellcolor{Gray1}0.91&\cellcolor{Gray1}0.47&&\cellcolor{Gray2}0.19&\cellcolor{Gray3}0.00&\cellcolor{Gray2}0.31\\
Rank 3&\cellcolor{Gray3}0.03&\cellcolor{Gray3}0.06&\cellcolor{Gray1}0.31&&\cellcolor{Gray1}0.91&\cellcolor{Gray1}0.81&\cellcolor{Gray2}0.19&&\cellcolor{Gray2}0.06&\cellcolor{Gray2}0.09&\cellcolor{Gray1}0.31&&\cellcolor{Gray4}0.00&\cellcolor{Gray4}0.00&\cellcolor{Gray2}0.19\\
Rank 4&\cellcolor{Gray1}0.94&\cellcolor{Gray1}0.88&\cellcolor{Gray2}0.16&&\cellcolor{Gray2}0.06&\cellcolor{Gray2}0.13&\cellcolor{Gray1}0.75&&\cellcolor{Gray3}0.00&\cellcolor{Gray3}0.00&\cellcolor{Gray3}0.09&&\cellcolor{Gray3}0.00&\cellcolor{Gray3}0.00&\cellcolor{Gray4}0.00\\
\hline
\end{tabular}
\end{table}

The results demonstrate that AILS-PR was ranked first in the vast  majority of the instances with regard to the average and best solutions. In X-n101-247, it presented the best  average solutions in 100\% of the instances and the top best solutions in 97\% of the instances. The top average and best solutions were achieved in, respectively, 97\% and 92\% in  X-n252-491; and, respectively, in 81\% and 100\% of the instances in  X-n502-1001. Regarding the average times, AILS-PR also outperformed the other methods with the best average times in 59\%, 44\% and 50\% of the respective group of instances X-n101-247, X-n252-491 and X-n502-1001.

In addition to the instances used in this experiment, Appendix~C shows the results of the literature solution methods and AILS-PR considering the instances proposed in \cite{Christofides1979} and \cite{Golden1998}.

\section{Final Remarks and Future Works}
\label{section:conclusao}

In this paper, we proposed a hybrid method to solve the Capacitated Vehicle Routing Problem (CVRP). First, we introduced a new version of the well-known ILS, here called Adaptive ILS (AILS). This method has an adaptive mechanism to control the diversity in the search process, by adjusting parameters related to the perturbation intensity and acceptance criterion. These components of ILS are highly important to the diversity control of the method.  The proposed approach can be used to solve other combinatorial optimization problems, since we propose AILS as a generic framework.

To propose a method that outperforms the state-of-the-art CVRP algorithms, we combined AILS with the widely employed intensification strategy known as Path-Relinking (PR). The proposed hybrid method was  named AILS-PR.  The hybridization guaranteed a good balance between diversification and intensification, since in comparison to AILS, it provided better and more robust results. 

In computational experiments, we tested  AILS-PR by employing a benchmark set containing 100 CVRP instances \citep{Uchoa2017}. A detailed analysis of the experiments was performed to state the better performance of AILS-PR in comparison to the UHGS, ILS-SP and SISRs. The results demonstrated that the proposed algorithm was significantly better than the state-of-the-art algorithms, being better in both solution quality and computational time. 

As future work, we suggest applying the proposed algorithm to other VRP problems, such as Multiple Depot VRP (MDVRP), Heterogeneous VRP (HVRP), Distance-Constrained VRP (DVRP) and VRP with Time Windows (VRPTW).

\section*{Acknowledgments}

The authors are grateful for the financial support provided by CNPq and FAPESP. Research carried out using the computational resources of the Center for Mathematical Sciences Applied to Industry (CeMEAI) funded by FAPESP (grants 2013/07375-0 and 2019/22067-6).

\bibliographystyle{elsarticle-harv} \bibliography{Bibliography}

\end{document}


\begin{frontmatter}

\title{A hybrid adaptive Iterated Local Search with diversification control to the Capacitated Vehicle Routing Problem\\(Online Appendices)} 

\author[label1]{Vinícius R. Máximo}
\ead{vinymax10@gmail.com}

\author[label1]{Mariá C. V. Nascimento\corref{cor1}}
\ead{mcv.nascimento@unifesp.br}
\cortext[cor1]{Corresponding author}

\address[label1]{Instituto de Ci\^{e}ncia e Tecnologia, Universidade Federal de S\~{a}o Paulo - UNIFESP\\ 
Av. Cesare G. Lattes, 1201, S\~{a}o Jos\'{e} dos Campos, SP, Brasil \\}

\end{frontmatter}

\appendix\section{AILS to solve the CVRP}

This section presents the introduced AILS to approach the CVRP. Each procedure of AILS for the CVRP is showed in the next sections.

\subsection{Construction of the Initial Solution}
\label{ConstrucaoInicial}

This section introduces the proposed heuristic to construct initial solutions  to the CVRP. In preliminary tests with AILS-PR,  we observed that the quality of the final solution of AILS-PR does not vary significantly considering different initial solutions. Therefore, as AILS-PR must be scalable, we suggest a simple heuristic to this step of the introduced algorithm, described in Algorithm \ref{construcaoinicial}. 

\begin{algorithm}[!htb]
\LinesNumbered
\SetAlgoLined
\KwData{Instance data}
\KwResult{A solution $s$}
$V_a\leftarrow V_c$\\
\For{$j=1,\ldots, \underline{m}$}{
 $s \leftarrow$ Construct $R^s_j=\{v_0^{j},v_1^{j}, v_0^{j}\}$ by randomly choosing $v_1^{j}$ from $V_a$\\
 $V_a \leftarrow V_a\backslash v_1^{j}$\\
} 
\While{$V_a\neq \emptyset$}
{ 
    Pick $v\in V_a$ at random and make     $V_a \leftarrow V_a\backslash v$\\
    Choose route $R^s_j$ where $v$ will belong considering $\mathscr{H}^a_1$ described in \ref{perturbacao}\\
    Add vertex $v$ in route $R^s_j$ at the lowest cost position\\
}
\Return{$s$}
\caption{Construct initial solution}
\label{construcaoinicial}
\end{algorithm}

According to Algorithm \ref{construcaoinicial}, the initial solution has $\underline{m}$ routes, each of them containing a single vertex chosen at random. The value of $\underline{m}$ is calculated according to Equation \eqref{minferior}, which is a lower bound on the number of routes.

 \begin{align}
&\underline{m}=\lceil \frac{1}{\bar{q}} \sum_{i \in V_c} q_i \rceil& \label{minferior}
\end{align}

Then, in random order, the other vertices are added to the routes. The insertion heuristic $\mathscr{H}^a_1$  described in \ref{perturbacao}  decides upon  which route to add each vertex as well as its position in the route. The resulting solution may not be feasible because the route can violate the capacity constraints when constructing the solution. For this reason, it is necessary to apply a feasibility algorithm  to the  solution afterwards, whose procedure will be discussed in  \ref{Factibilizacao}.

\subsection{Perturbation Strategy}
\label{perturbacao}

The ILS perturbation process consists of producing changes in a solution $s$ by a heuristic to obtain a neighbor solution $s'$ different from $s$. It has two steps. The first, which we call removal heuristics, consists of deciding which vertices are to be removed from the current solution. The second step, called insertion heuristics, defines where the vertices removed in the first step are to be inserted into the solution. In this paper, we adopted three removal heuristics and two insertion heuristics. The perturbation methods adopted in the proposed AILS-PR do not guarantee the feasibility of the resulting solution.  

\textbf{Concentric removal, $\mathscr{H}^r_1$: } In this heuristic, a randomly selected  vertex $v_r$ is removed from the solution as well as the $\omega_{\mathscr{H}_1}-1$  closest vertices to $v_r$. In other words, this heuristic removes a set of vertices belonging to a coverage radius whose center is represented by vertex $v_r$.

\textbf{Proximity removal, $\mathscr{H}^r_2$: } This heuristic consists of removing $\omega_{\mathscr{H}_2}$ vertices from the solution. To select which vertices to remove, the strategy uses the concept of a proximity index. The proximity index is a measure that indicates how close a particular vertex is to its route. To estimate the proximity,  we do not consider the Euclidean distance metric, because it has an intrinsic magnitude of the distribution of points in the space. Instead,  we use a relative proximity metric, which allows a comparison among all vertices independent of the magnitude of the pairwise distance between points.

Let $\Pi_R: V_c \rightarrow \{0,\ldots,n-1\}^{|R|-2}$ be the function that assigns to a given vertex $v$ the order of the vertices in a given route $R\backslash\{0,v\}$ according to their Euclidean distance to $v$. For example, let $R=\{0,1,4,6,7,0\}$ be a route of a graph with 10 vertices, i.e, $V=\{0,1,\ldots, 9\}$. The length of $R$, here denoted by $|R|$, is $5$. To calculate $\Pi_R(4)$, one must consider the distances between 4 and all vertices in the route. Suppose that vertices 1, 6 and 7 are ranked as 3, 5 and 8 closest vertices to 4 among all vertices in the graph, respectively. Therefore, $\Pi_R(4)=\{3,5,8\}$. In addition, consider $minset(a,S)$ the sum of the $a$ smallest elements of a set $S$. Thus, the proximity index of a vertex $v$ with respect to route $R^j_s$ where it belongs to in a given solution $s$ is: 

\begin{equation}prox(v,\rho,R^s_j) = \frac{minset(\min\{\rho,|R^s_j|-2\},\Pi_{R^s_j}(v))}{\min\{\rho,|R^s_j|-2\}}\label{proximidade} \end{equation}

In this proximity index, we consider only the $ \rho $  closest vertices to $ v $. 
An integer random value for $ \rho $ is chosen randomly at each iteration of the algorithm in the interval $ [1, \lfloor n / m \rfloor] $. In addition, it is possible to observe that $prox(v,\rho,R^s_j) \in [1,n-1]$.

The proposed vertex removal strategy  prioritizes the removal of vertices with the highest proximity indexes. To accomplish this task, we define $S $ as the set of all vertices belonging to the solution which have not  been removed yet. Next, we sort all the elements of $ S $ in decreasing order of proximity index  and remove $ \omega_{\mathscr{H}^r_2} $ vertices from $ S $ according to the probability $ P (v) = \frac{2 (| S | -o_v) -1} {| S |^2} $, where $o_v$ represents the vertex position $ v $ in the sorting of $S$. This probability distribution returns the probability of removing a vertex $ v $ which belongs to $S$. The lower the position  $ o_v \in [0, | S | -1] $ the greater the probability of removing $ v $.

\textbf{Removal of Vertex Sequences, $\mathscr{H}^r_3$: } In this removal strategy, some vertex sequences belonging to a route are removed. The sequences can have different sizes and are chosen randomly from the sequences of vertices belonging to the solution. The number of vertices to be removed is controlled by the variable $ \omega_{\mathscr{H}^r_3} $ and the length of the sequence varies randomly in the interval $ [1, \lfloor n / m\rfloor] $ since the total number of vertices already removed does not exceed $ \omega_{\mathscr {H}^r_3} $. Sequences can contain the depot, however, it will not be removed from the solution. For example, in route $ R = \{0,2,6,8,11,0 \} $, a possible sequence would be $ (8,11,0,2) $. In this case, only three elements could be removed by this heuristic: 8, 11 and 2.

The two insertion heuristics employed in this study are presented next.

\textbf{Insertion by Proximity, $\mathscr{H}^a_1$:} This heuristic consists of inserting a vertex $ v $ in the lowest proximity index route  calculated according to Equation \eqref{proximidade}. Accordingly, vertex $ v $ will be inserted in route $ R^s_{\hat{j}} = \arg \min_{R ^ s_j \in \mathscr{R}} prox (v, \rho, R^s_j ) $, in the lowest cost position  $ \hat{i} $ according to Equation~\eqref{bestPosition}. The insertion position of $ v $ will be between the two vertices $ v_i^{\hat{j}} $ and $ v_{i + 1}^{\hat{j}} $ that results in the least insertion cost in the route $ R^s_{\hat{j}} = \{ v_0^{\hat{j}}, v_1^{\hat{j}}, \ldots, v^{\hat{j}}_{m_{\hat{j}}} \} $.

\begin{align}
\hat{i}=\arg\min_{i\in\{0,1,\ldots,m_{\hat{j}}\}} d(v_i^{\hat{j}},v)+d(v_{i+1}^{\hat{j}},v)-d(v_i^{\hat{j}},v_{i+1}^{\hat{j}}) & \label{bestPosition}
\end{align}

\textbf{Insertion by Cost, $ \mathscr{H}^a_2 $:} In this heuristic,  vertices are inserted in the routes that have the lowest insertion cost, that is, vertex $ v $ will be inserted in the route $ R^s_{\hat{j}} = \arg \min_{R^s_j \in \mathscr{R}} c (v^j_{\hat{i}}, v^j_{\hat{i} +1} , v) $. The position of insertion of  vertex $ v $ in the route $ R^s_{\hat {j}} $ follows the criterion of lowest cost.

\begin{algorithm}[!htb]
\LinesNumbered
\SetAlgoLined
\KwData{Solution $s$ and $\mathscr{H}^r_k$}
\KwResult{Solution $s'$}
$s' \leftarrow s$\\
Vary the number of routes $m$ of solution  $s'$ in 1 unity with probability $1/\gamma$ as long as  $m \geq \underline{m}$\\
Choose  an integer number at random within interval $ [1,\lfloor n/m\rfloor]$ and assign to  $\rho$ \\
Remove $\omega_{\mathscr{H}^r_k}$ vertices from solution $s'$ according to heuristic \textbf{ $\mathscr{H}^r_k$}\\
 Randomly selected one of the insertion heuristics $\{\mathscr{H}^a_1,\mathscr{H}^a_2\}$ to add to $s'$ all vertices removed\\
\caption{Perturbation Procedure}
\label{pertubacao}
\end{algorithm}

In Algorithm ~\ref{pertubacao}, after $s$ is copied to $s^{'}$, the number of routes of solution $s'$ is defined. The number of routes can vary by one unit from the input solution, $s$, with probability $1/\gamma$ -- being $\gamma$ a parameter of AILS-PR -- as long as the resulting $m$ is greater  than or equal to $\underline{m}$. If from the number of routes $m$ is subtracted  1, a randomly chosen route must be removed from solution $s'$. If $m$ is incremented by 1, a new empty route will be inserted into the solution. After this step, we choose an integer random value for $ \rho \in [1, \lfloor n / m \rfloor] $ and start the process of removing and inserting vertices in the routes. In this way, the strategy removes $ \omega_{\mathscr{H}^r_k} $ vertices from solution $ s' $  according to heuristic \textbf{$ \mathscr{H}^r_k$}. Then, the vertices removed from the solution are inserted  by an insertion heuristic.

\subsection{Feasibility and Local Search Strategies}

In this section, the feasibility and local search strategies are thoroughly explained.  Both feasibility and local search strategies perform moves of vertices between routes.  They use the neighborhood structure described next.


\noindent\textbf{Neighborhood structure of moves between routes}

The inter-route neighborhood investigated of the strategies is formed by the movements of $ \lambda $-interchange proposed by \citet{Osman1993} and 2-opt$^*$ \citep{Potvin1995}. 
The value of $ \lambda $ establishes the maximum size of adjacent vertices to be moved. After parameter set up  experiments, we adopt $ \lambda = 1 $. Therefore, the neighborhood structure will be composed of movements that swap positions of two vertices -- called $ \mathbf{Swap} $ -- and  change the position of a vertex -- known as $ \mathbf{Shift} $.

Let the set of moves $\mathscr{N}= \{\mathscr{N}_1, \mathscr{N}_2, \mathscr{N}_3\}$ define the neighborhood structure of the proposed strategy. Each move is described as follows. 

\noindent\textbf{$\mathbf{Shift}, \mathscr{N}_1$}: In this move, given a solution $s$, a neighbor $s'$ is obtained by moving a vertex $v_k^i$ from a certain route $R_{i}^s$ to another route $R_{j}^s$, placing it between vertices $v^j_l$ and $v^{j}_{l+1}$ where $i\neq j$, $k>0$ and $l<|R^s_j|$. The objective function value of $s'$ can be derived by evaluating the variation in the objective function value of $s$:

\begin{equation}
  \Delta^{\mathscr{N}_1}_{v_k^i,v^j_l}(s,s') =
        d(v^i_{k-1},v^i_{k+1})+d(v^{j}_l,v_k^i)+d(v_k^i,v^{j}_{l+1})-d(v^i_{k-1},v_k^i)-d(v_k^i,v^i_{k+1})-d(v^{j}_l,v^{j}_{l+1}) 
 \end{equation}


If $ \Delta^{\mathscr{N}_1}_{v_k^i,v^j_l}(s,s')$ is negative, it means that $s'$ -- which is a neighbor of $s$ -- is a better solution than $s$.

\textbf{$\mathbf{Swap}, \mathscr{N}_2$}: In this move,  vertices are swapped in a pair of routes.  Given a solution $s$, a neighbor $s'$ is  obtained by swapping positions between  vertices $v_k^i$  from route $R_{i}^s$ with vertex $v^{j}_{l}$ from route $R_{j}^s$. In this case,  $i\neq j$,  $k>0$ and $l<|R_{j}^s|$. The objective function value of $s'$ after this move can be derived by evaluating the variation in the objective function value of $s$:

\begin{equation}
    \begin{split}
  \Delta^{\mathscr{N}_2}_{v_k^i,v^j_l}(s,s') = 
        d(v^i_{k-1},v^j_{l})+d(v^{j}_{l},v^i_{k+1})+d(v^{j}_{l-1},v^i_k)+d(v^i_k,v^{j}_{l+1})\\
        -d(v^i_{k-1},v_k^i)-d(v_k^i,v_{k+1}^i)-d(v^{j}_{l-1},v^{j}_{l})-d(v^{j}_{l},v^{j}_{l+1})
    \end{split}
\end{equation}

\textbf{$\mathbf{2-opt^*}, \mathscr{N}_3$: } This move consists of exchanging two pieces of routes. It has as reference two edges, one from a different route, which will be exchanged to define a neighbor solution. Let $R_{i}^s$ and $R_j^s$ be the routes of solution $s$ whose edges $(v^i_k,v_{k+1}^i)$ and $(v^j_l,v^j_{l+1})$, respectively, will be exchanged to produce solution $s'$.
 
 \begin{equation}
  \Delta^{\mathscr{N}_3}_{v_k^i,v^j_l}(s,s') = 
        d(v^{j}_{l},v_{k+1}^i)+d(v^i_k,v^{j}_{l+1})-d(v_{k}^i,v^i_{k+1})-d(v^{j}_{l},v^{j}_{l+1})   
 \end{equation}

In our implementation, all moves are made in $O(1)$.

\noindent\textbf{Neighborhood structure of intra-route moves}

In addition to the inter-route moves introduced in the previous section, in AILS, another neighborhood structure defines a strategy called ``intra-local search''. The intra-route moves employed in the intra-local search consist of promoting changes in the position of the vertices that belong to the same route. The set of moves applied in a route allows to establish the neighborhood $\mathscr{N}^- = \{\mathscr{N}_1^-, \mathscr{N}_2^-, \mathscr{N}_3^-\}$ composed by 3 moves:

$\mathbf{Shift}^{-}, \mathscr{N}_1^-$: A vertex is relocated to another position of the same route.

$\mathbf{Swap}^{-}, \mathscr{N}_2^-:$ This move swaps two vertices in the same route.

$\mathbf{2-opt}^ {-}, \mathscr{N}_3^-:$ This move removes a pair of non-adjacent edges and adds two new edges in the route to guarantee the feasibility of the route.

It is possible to calculate the variation in the solution values to identify the neighboring solution locally similar to the movements between routes just explained.

\noindent\textbf{Intra-Route Local Search Heuristic}
\label{BLIntra}

The intra-route local search used in the AILS-PR is performed following the best improvement strategy according to the neighborhood $\mathscr{N}^-$.  The algorithm searches for the best neighbor $s' \in \mathscr{N}^-$ found by applying to $s$ the movement in $R^{s}_i$. If this neighbor is found and  $\Delta^{\mathscr{N}^-}_{v_k^i,v^i_l}(s,s')<0$, solution $s$ is updated. The search process repeats until $s$ becomes a local optimum. It is worth pointing out that the strategy  accepts only movements whose costs are lower than zero, which is,  accepts only better quality solutions than the current one.


The feasibility  and local search procedures have the same main structure. They  differ in some conditions to be satisfied to a move to be performed. In line with this, the feasibility heuristic aims at finding a feasible good quality solution. The goal of the local search procedure is to reach the local optimum.  The neighborhood selection strategy follows the best improvement strategy. Algorithm~\ref{alg:neighborhoodsearch} shows the neighborhood search algorithm that, if the input solution $s$ is infeasible, will work as a feasibility heuristic, otherwise, it will be a local search. 



\begin{algorithm}[!htb]
\LinesNumbered
\SetAlgoLined
\KwData{Solution $s$, $feasible$}
\KwResult{Solution $s$}
$LM \leftarrow \emptyset$\\ 
\Repeat{$s$ is feasible}
{
    \ForEach{ $R^{s}_i \in \mathscr{R}^s $}
    {
        Update $LM$ ($R^{s}_i$, $LM$, $s$, $feasible$)\\
    }
   \While{$|LM|>0$}
   {
        Find $\Theta^{\mathscr{N}}_{v^i,v^j}(s,s')\in LM$ \label{linha7}\\
        $s \leftarrow s'$\\
        $LM \leftarrow LM \setminus \{\Theta^{\mathscr{N}}_{v^i,v^j}(s,x) \in  LM\}$\\
        Apply the intra-route local search to routes $R^{s}_i$ and $R^{s}_j$\\
        Update $LM$  ($R^{s}_i$, $LM$, $s$, $feasible$) and Update $LM$  ($R^{s}_j$, $LM$, $s$, $feasible$)\\
    }
    \If{$s$ is infeasible}
    {
        Add a new route $\mathscr{R}^s \leftarrow \mathscr{R}^s \cup R^s_{m+1}$\\
    }
}
\caption{Neighborhood Search}
\label{alg:neighborhoodsearch}
\end{algorithm}

In Algorithm~\ref{alg:neighborhoodsearch}, the input data are both the solution $s$ and a boolean variable that indicates whether or not such a solution is feasible. First, the list of movements, referred to as  $ LM $,  is initialized as an empty set. The main loop occurs once for the local search, since the input solution is feasible, and while $s$ is infeasible  in the feasibility heuristic. Then, the method performs a search in the neighborhood of $s$ as described in Algorithm \ref{alg:LM}. This algorithm returns the set $LM$ which contains all the movements between  vertices $ v_k^i \in R^s_i $ and $ v_l^j \in \delta (v_k^i) ~ | ~ i \neq j $, where $ \delta (v_k^i) $ represents the set of closest vertices  to $ v_k^i$ according to the proximity index. The size of this set is limited by parameter $ \varphi \in [1, n-1] $, which means that $ | \delta (v) | = \varphi, \forall ~ v \in V $. After this first step of verifying the movements, the algorithm updates the current solution $s$ while there are movements to perform in list $ LM $. In this step, the algorithm updates the current solution in the best neighbor contained in $LM$. This neighbor is represented by the movement $ \Theta^{\mathscr{N}}_{v^i, v^j} (s, s') $. Then,  $ LM $ is updated, removing all movements involving the routes $ R^{s}_i $ and $ R^{s}_j $, as there have been changes in these routes. Then, an intra-local search is applied to $s$ and   list $ LM $ is updated.
If  $ LM $  is empty and  $ s $  remains infeasible, then a new route is added to the solution and the loop is repeated.

\begin{algorithm}[!htb]
\LinesNumbered
\SetAlgoLined
\KwData{ $R^{s}_i$, $LM$, solution $s$, $feasible$}
\KwResult{Updated $LM$}
\ForEach{$v_k^i \in R^s_i$}
{
    \ForEach{$v_l^j \in \delta(v_k^i) ~|~ i \neq j$ and $R^{s}_i$ and $R^{s}_j$ meet feasibility condition \label{linha2} }
    {
        \ForEach{ $\mathscr{N}_k \in \mathscr{N}$}
        {
     \If{$\mathscr{N}_k$ provides a solution $s''$ that respects feasibility/improvement conditions}{    \label{linha4}        
     \If{ $  \exists~ \Theta^{\mathscr{N}}_{v^i, v^j} (s, s') \in LM$ but $s''$ is better }
            {\label{linha5}
                
                    $\Theta^{\mathscr{N}}_{v^i,v^j}(s,s') \leftarrow \Theta^{\mathscr{N}_k}_{v_k^i,v^j_l}(s,s'')$\\
                }
                
                \If{$\not\exists~ \Theta^{\mathscr{N}}_{v^i,v^j}(s,s') \in LM $ }
                {
                   $LM \leftarrow LM \cup \Theta^{\mathscr{N}_k}_{v_k^i,v^j_l}(s,s'')$\\
                }
            }
        }
    }
}
\caption{Update $LM$}
\label{alg:LM}
\end{algorithm}

Next sections present the conditions and criteria employed to the feasibility and local search procedures, respectively.



\subsubsection{Criteria for the feasibility procedure}
\label{Factibilizacao}

In the CVRP, the capacity constraints are the only set of constraints that can be violated in the search process of the introduced algorithm. Therefore, the feasibility heuristic performs vertex movements -- neighborhood  $ \mathscr{N} $ -- between routes to adjust the routes that exceed the limit capacity. 

Let $\mathscr{R}=\{R^{s}_1,R^{s}_2, \ldots,R^{s}_m\}$ be the routes of a solution $s$. Consider  $slack(R^{s}_j)$ the  slack in the load constraints  of route $R^{s}_j$ which can be calculated by Equation~\eqref{slack}. $R^{s}_j$ is feasible iff $slack(R^{s}_j)\geq 0$.

\begin{align}
&slack(R^{s}_j)=\bar{q}-\sum_{i \in R^{s}_j} q_i & \label{slack}
\end{align}

The moves of vertices/edges between routes $R^{s}_i,R^{s}_j  \in \mathscr{R}$ where $R^{s}_i$ is infeasible and $R^{s}_j$ is feasible can provide a neighbor solution $s' \in \mathscr{N}(s)$ with reduction in the violation of the load constraints. The reduction occurs in two particular cases:

\begin{enumerate}
\item  When the slack in the load constraints in  route  $R^{s'}_i$ is greater than the slack of route $R^{s}_i$ and the route  $R^{s'}_j$ is feasible. In other words, when $slack(R^{s'}_i) > slack(R^{s}_i)$  and $slack(R^{s'}_j)\geq 0$. 
\item When route $R^{s'}_i$ becomes feasible and the slack in the load constraints in  $R^{s'}_j$ is greater than the slack of $R^{s}_i$. In other words, when  $slack(R^{s'}_i)\geq 0$ and $ slack(R^{s'}_j) > slack(R^{s}_i)$.  
\end{enumerate}

Let $\Omega(s,s')$ be the ``feasibility gain'' after a move applied to solution $s$ to obtain neighbor $s' \in \mathscr{N}(s)$. The value of  $\Omega(s,s')$ is calculated according to Equation~\eqref{ganho}. Therefore, if $\Omega(s,s') >0$, the move provided a $s'$ with load constraints violation lower than the load constraints violation in $s$.

\begin{align}
\scalefont{0.8}
\Omega(s,s') = \min\{0, slack(R^{s'}_i)\} + \min\{0, slack(R^{s'}_j)\} - \min\{0, slack(R^{s}_i)\} - \min\{0, slack(R^{s}_j)\}
& \label{ganho}
\end{align}

The decision on which moves to perform follows a pair of criteria: i) the cost of the movement -- this criterion is considered when the cost of the movement is less than or equal to 0; ii) the cost of the movement and the feasibility gain -- when the cost of the movement is greater than 0.  To meet these two criteria, we evaluate the value of $ \Lambda^{\mathscr{N}_k}_{v_k^i,v^j_l}(s,s') $ considering a given move $\mathscr{N}_k$ from route $R_{i}^s$ to route $R_{j}^s$ according to Equation~\eqref{custoFac}. Thus, the lower the value of $ \Lambda^{\mathscr{N}_k}_{v_k^i,v^j_l}(s,s')$, the better  is the move.

\begin{align}
\Lambda^{\mathscr{N}_k}_{v_k^i,v^j_l}(s,s') = \left \{ \begin{array}{ll}
\Delta^{\mathscr{N}_k}_{v_k^i,v^j_l}(s,s'), & \textrm{if  $\Delta^{\mathscr{N}_k}_{v_k^i,v^j_l}(s,s') \leq 0$}  \\\\
\frac{\Delta^{\mathscr{N}_k}_{v_k^i,v^j_l}(s,s')}{\Omega(s,s')}, & \textrm{if $\Delta^{\mathscr{N}_k}_{v_k^i,v^j_l}(s,s') > 0$} 
\end{array}\right. & \label{custoFac}
\end{align}

In summary, the particularities of the feasibility procedure in Algorithms~\ref{alg:neighborhoodsearch} and \ref{alg:LM} are:

\begin{itemize}
    \item In line \ref{linha7} of Algorithm~\ref{alg:neighborhoodsearch},  the move $\Theta^{\mathscr{N}}_{v^i,v^j}(s,s')$ chosen is $ \arg \min_{\Theta(s,s'') \in LM} \Lambda(s,s'')$;
    \item The feasibility condition of routes $R^{s}_i$ and $R^{s}_j$ in line~\ref{linha2} of Algorithm~\ref{alg:LM} is: $(slack(R^{s}_i)\geq0~XOR~slack(R^{s}_j)\geq0)$;
    \item Condition in line~\ref{linha4} of Algorithm~\ref{alg:LM} is $\Omega^{\mathscr{N}_k}_{v_k^i,v^j_l}(s,s'') > 0$;
    \item Condition in line~\ref{linha5} of Algorithm~\ref{alg:LM} is $\exists~ \Theta^{\mathscr{N}}_{v^i,v^j}(s,s') \in LM ~|~ \Lambda^{\mathscr{N}_k}_{v_k^i,v^j_l}(s,s'') < \Lambda^{\mathscr{N}}_{v^i,v^j}(s,s') $.
\end{itemize}

The next section presents the criteria employed for the local search of the proposed method.

\subsection{Criteria for the local search procedure}
\label{buscaLocal}

The difference between the feasibility and  local search criteria of AILS-PR resides in the chosen movements. In the local search are evaluated considering every pair of routes since they are all viable. In addition, in local search, only the costs of the optimization process are considered. 

In line with this, the particularities of the local search procedure in Algorithms~\ref{alg:neighborhoodsearch} and \ref{alg:LM} are:

\begin{itemize}
    \item In line \ref{linha7} of Algorithm~\ref{alg:neighborhoodsearch},  the move $\Theta^{\mathscr{N}}_{v^i,v^j}(s,s')$ chosen is  $\arg \min_{\Theta(s,s'') \in LM} \Delta(s,s'')$;
    \item There is no feasibility condition of routes $R^{s}_i$ and $R^{s}_j$ in line~\ref{linha2} of Algorithm~\ref{alg:LM}, so, the search considers every pair of routes;
    \item Condition in line~\ref{linha4} of Algorithm~\ref{alg:LM} is $\Omega^{\mathscr{N}_k}_{v_k^i,v^j_l}(s,s'') \geq 0$ and $\Delta^{\mathscr{N}_k}_{v_k^i,v^j_l}(s,s'') < 0$;
    \item Condition in line~\ref{linha5} of Algorithm~\ref{alg:LM} is $\exists~ \Theta^{\mathscr{N}}_{v^i,v^j}(s,s') \in LM ~|~ \Delta^{\mathscr{N}_k}_{v_k^i,v^j_l}(s,s'') < \Delta^{\mathscr{N}}_{v^i,v^j}(s,s') $.
\end{itemize}


 Both the feasibility heuristic and the local search algorithm have $ O(n \varphi) $ complexity. To sum up, we present in Algorithm~\ref{AILSCVRP} the AILS proposed to solve the CVRP.

\begin{algorithm}[!htb]
\LinesNumbered
\SetAlgoLined
\KwData{Instance data}
\KwResult{The best solution found $s^*$}
$s^r \leftarrow$ Construct initial solution -- Algorithm \ref{construcaoinicial}\\ 
$s^r,s^* \leftarrow$ Neighborhood Search($s^r$,0) -- Algorithm~\ref{alg:neighborhoodsearch}\\ 
$it\leftarrow 1$\\
\Repeat{$it\geq it_{max}$}
{ 
    $\mathscr{H}^r \leftarrow$ Choose one removal heuristic $\mathscr{H}^r \in \{\mathscr{H}^r_1,\mathscr{H}^r_2,\mathscr{H}^r_3\}$, explained in \ref{perturbacao} \\
$s \leftarrow$ Perturbation Procedure($s^r$, $\mathscr{H}^r$) --- Algorithm~\ref{pertubacao}\\
	   $s \leftarrow$ Neighborhood Search($s$,0)  -- Algorithm~\ref{alg:neighborhoodsearch}\\
    $s \leftarrow$ Neighborhood Search($s$,1)  -- Algorithm~\ref{alg:neighborhoodsearch}\\
	Update the diversity control parameter $\omega_{\mathscr{H}^r_k}$ considering the distance between $s$ and $s^r$\\
    $s^r\leftarrow $ Apply acceptation criterion to $s$\\
	Update the acceptance criterion\\
	Assign $s$ to $s^*$ if $ f(s) < f(s^*)$\\
    $it \leftarrow it+1$\\
}
\caption{AILS}
\label{AILSCVRP}
\end{algorithm}

\section{Calculating Priority \texorpdfstring{$p_v$}{Lg} to perform PR}
\label{prioridadePv}


The priority $ p_v $ represents the priority of  vertex $ v \in NF $ to be moved to another route. The calculation of the priority depends on the criterion $ C \in \mathscr{C} $  and the states of $v$ considering its origin and destination routes, respectively,  $ R^{s_i}_{o} $ and  $ R^{s_i}_{d}$.
The states are determined according to the feasibility of these origin and destination routes. 
Let  $ R^{s_i}_{o}(v^+) $ and $ R^{s_i}_{o} (v^-) $ be the route $ R^{s_i}_{o} $ before and after removing $ v $, respectively. Consider also  $ R^{s_i}_{d} (v^-) $  and $ R^{s_i}_{d} (v^+) $ the route $ R^{s_i}_{d} $ before and after receiving  vertex $ v $ a route $ R^{s_i}_{d} $, respectively. 
Both  $ R^{s_i}_{o} $  and  $ R^{s_i}_{d} $  have 3 possible states regarding feasibility. These states are presented in Table \ref{6estados} which also indicates how to calculate  the priority of vertex $v$.

\begin{table}[!htb]
\caption{List with the six states used to calculate priority $p_v$.}
\label{6estados}
\begin{tabular}{|c|c|c|c|c|c|c|c|}
\hline
\multicolumn{4}{|c|}{Origin route ($O$)}  &\multicolumn{4}{|c|}{Destination route ($D$)} \\ \hline
State  & \multicolumn{1}{|c|}{$R^{s_i}_{o}(v^+)$}  & \multicolumn{1}{|c|}{$R^{s_i}_{o}(v^-)$}   & Priority factor & \multicolumn{1}{|c|}{State}  & \multicolumn{1}{|c|}{$R^{s_i}_d(v^-)$} &  \multicolumn{1}{|c|}{$R^{s_i}_d(v^+)$} & Priority factor\\
\hline
\rowcolor{Gray}
1 & Infeasible & Feasible      & $+ 1$  & 4 & Feasible  & Infeasible & $- 1$  \\
2 & Infeasible & Infeasible    & $+ 1$ & 5 & Infeasible & Infeasible &  $ - 1$ \\
\rowcolor{Gray}
3 &  Feasible & Feasible & $- 1$ & 6 & Feasible & Feasible &  $ + 1$ \\
\hline
\end{tabular}
\end{table}

States 1, 2 and 6 have an increased priority because we consider that these conditions are favorable to find a feasible solution by PR. States 3, 4 and 5 have their priority decreased.

Let $\mathscr{C}=\{\mathscr{C}_1, \ldots, \mathscr{C}_{10}\}$ be the set of criteria used in PR. Each criterion represents a specific combination of active states. In this way, the priority $p_v$ can be calculated in different ways, thus allowing the PR to go through different paths between the solutions. Considering all possible combinations of the six states, there are $2^6=64$ possible different criteria. However, after computational experiments, we chose to use only the 10 criteria presented in Table \ref{criterios}. 

\begin{table}[]
\centering
\caption{Description of the states employed in each criterion.}
\label{criterios}
\begin{tabular}{|c|c|c|c|c|c|c|c|c|c|c|}
\hline
\multirow{2}{*}{State}   & \multicolumn{10}{c|}{Criterion}\\
\cline{2-11}
&$\mathscr{C}_1$ &$\mathscr{C}_2$ &$\mathscr{C}_3$ &$\mathscr{C}_4$ &$\mathscr{C}_5$ &$\mathscr{C}_6$ &$\mathscr{C}_7$ &$\mathscr{C}_8$ &$\mathscr{C}_9$ &$\mathscr{C}_{10}$ \\
\hline
\rowcolor{Gray}
1& x &   & x & x & x & x &   &   &   &   \\
2&   &   &   & x & x & x &   & x & x &  \\
\rowcolor{Gray}
3& x & x & x & x &   & x & x & x & x & x \\
4&   &   &   &   &   &   &   &   &   &x   \\
\rowcolor{Gray}
5& x &   &   & x & x &   & x & x &  &  \\
6& & & & & & & & & & x\\
\hline
\end{tabular}
\end{table}

We define $O =\{1, 2, 3 \} $ and $D=\{4, 5, 6 \}$ as being the sets containing the states of the origin and destination routes, respectively. Therefore, to calculate the priority $p_v$, it is necessary to consider the state pair $(o, d)$, where $o \in O$ and $d \in D$. The priority calculation as indicated in Table \ref{6estados} will be performed if, and only if, the state is active in the chosen criterion. Each priority value of a move will be the result of 0, 1, or 2 operations of decrement/increment. Considering the available updating possibilities, the priority may have the following values: $\{- 2, -1, 0, 1, 2 \}$.
Table \ref{prioridade} presents the priority values for each  state pair $(o, d)$. These values were derived considering the operations contained in Table~2, starting the priority values at zero.

\begin{table}[!ht]
\centering
\scalefont{0.7}
\caption{Priority value of each criterion considering every combination of states.}
\label{prioridade}
\begin{tabular}{|c|c|c|c|c|c|c|c|c|c|c|}
\hline
\multirow{2}{*}{Pair of states} &\multicolumn{10}{c|}{Criterion} \\ \cline{2-11}

&$\mathscr{C}_1$ &$\mathscr{C}_2$ &$\mathscr{C}_3$ &$\mathscr{C}_4$ &$\mathscr{C}_5$ &$\mathscr{C}_6$ &$\mathscr{C}_7$ &$\mathscr{C}_8$ &$\mathscr{C}_9$ &$\mathscr{C}_{10}$ \\
\rowcolor{Gray}
$(1,4)$ &1&0&1&1&1&1&0&0&0&0\\
$(1,5)$ &0&0&1&0&0&1&-1&-1&0&-1 \\
\rowcolor{Gray}
$(1,6)$ &1&0&1&1&1&1&0&0&0&1\\
$(2,4)$ &0&0&0&1&1&1&0&1&1&0\\
\rowcolor{Gray}
$(2,5)$ &-1&0&0&0&0&1&-1&0&1&-1\\
$(2,6)$ &0&0&0&1&1&1&0&1&1&1\\
\rowcolor{Gray}
$(3,4)$ &-1&-1&-1&-1&0&-1&-1&-1&-1&-1\\
$(3,5)$ &-2&-1&-1&-2&-1&-1&-2&-2&-1&-2\\
\rowcolor{Gray}
$(3,6)$&-1&-1&-1&-1&0&-1&-1&-1&-1&0\\
\hline
\end{tabular}
\end{table}

\section{Additional Experiment}

This appendix presents an experiment carried out with AILS-PR to solve the instances in \citep{Christofides1979}. 

The values reported in the Table~\ref{resultV} are:

\begin{itemize}
    \item BKS: the objective function value of the best-known solutions. The values highlighted with $*$ are optimal solutions according to the information available in the repository \url{http://vrp.atd-lab.inf.puc-rio.br}\footnote{Most recent data as of 14 november 2020 - 18.00GMT.}.
    
    \item Avg: reports the average objective function value of the solutions obtained in  50 rounds.
    
    \item \textit{gap}: represents the mean gap of the values of the Avg column relative to the BKS and is calculated according to Equation \eqref{gap}. 
    \begin{equation}\label{gap}
      gap=100(f(s)-\mbox{BKS})/\mbox{BKS}
\end{equation}

    \item Best: objective function value of the best solution obtained in  50 runs.
    
    \item Time: represents the average time in minutes that the algorithm took to find the best solution in each run.
\end{itemize}

The background of the top results for each instance is highlighted in dark gray; the second best results in medium gray; and the third best results in light gray. The configuration of the computer where the experiments of the reference algorithms were run are provided in the end of Table \ref{resultV}. The results of UHGS \citep{Vidal2014} are not reported in the literature, the reason why they do not appear in the table.

\addtolength{\tabcolsep}{-3pt}  \begin{table}[!ht]
\centering
\scalefont{0.6}
\caption{Results of experiments with instances proposed in \citep{Christofides1979} and \cite{Golden1998}.}
\label{resultV}
\centering
\begin{threeparttable}
\begin{tabular}{ccccccccccccccccc}
\hline
&& \multicolumn{3}{c}{ILS-SP \citep{Subramanian2013}} && \multicolumn{3}{c}{SISRs \citep{Christiaens2020}} && \multicolumn{3}{c}{AILS-PR}\\ 
\cline{3-5} \cline{7-9} \cline{11-13} 
Instance&BKS&Avg~(\textit{gap})&Best&Time\tnote{1} &&Avg&Best&Time\tnote{2}&&Avg~(\textit{gap})&Best&Time\tnote{3} \\
\hline
CMT1&524.61&\cellcolor{Gray1}524.61~(0.0000)&\cellcolor{Gray1}524.61&\cellcolor{Gray3}1.48&&\cellcolor{Gray1}524.61~(0.0000)&\cellcolor{Gray1}524.61&\cellcolor{Gray1}0.00&&\cellcolor{Gray1}524.61~(0.0000)&\cellcolor{Gray1}524.61&\cellcolor{Gray2}0.30\\
CMT2&835.26&\cellcolor{Gray1}835.26~(0.0000)&\cellcolor{Gray1}835.26&\cellcolor{Gray2}13.52&&\cellcolor{Gray3}835.4~(0.0168)&\cellcolor{Gray1}835.26&\cellcolor{Gray3}57.60&&\cellcolor{Gray1}835.26~(0.0000)&\cellcolor{Gray1}835.26&\cellcolor{Gray1}4.70\\
CMT3&826.14&\cellcolor{Gray1}826.14~(0.0000)&\cellcolor{Gray1}826.14&\cellcolor{Gray1}12.49&&\cellcolor{Gray1}826.14~(0.0000)&\cellcolor{Gray1}826.14&\cellcolor{Gray3}69.60&&\cellcolor{Gray1}826.14~(0.0000)&\cellcolor{Gray1}826.14&\cellcolor{Gray2}24.86\\
CMT12&819.56&\cellcolor{Gray1}819.56~(0.0000)&\cellcolor{Gray1}819.56&\cellcolor{Gray3}5.23&&\cellcolor{Gray1}819.56~(0.0000)&\cellcolor{Gray1}819.56&\cellcolor{Gray2}4.80&&\cellcolor{Gray1}819.56~(0.0000)&\cellcolor{Gray1}819.56&\cellcolor{Gray1}0.03\\
CMT11&1042.12&\cellcolor{Gray1}1042.12~(0.0000)&\cellcolor{Gray1}1042.12&\cellcolor{Gray2}20.66&&\cellcolor{Gray1}1042.12~(0.0000)&\cellcolor{Gray1}1042.12&\cellcolor{Gray3}108.00&&\cellcolor{Gray1}1042.12~(0.0000)&\cellcolor{Gray1}1042.12&\cellcolor{Gray1}0.58\\
CMT4&1028.42&\cellcolor{Gray2}1028.73~(0.0301)&\cellcolor{Gray1}1028.42&\cellcolor{Gray2}53.48&&\cellcolor{Gray3}1030.68~(0.2198)&\cellcolor{Gray1}1028.42&\cellcolor{Gray3}427.80&&\cellcolor{Gray1}1028.42~(0.0000)&\cellcolor{Gray1}1028.42&\cellcolor{Gray1}46.01\\
CMT5&1291.29&\cellcolor{Gray3}1293.18~(0.1464)&\cellcolor{Gray1}1291.45&\cellcolor{Gray2}625.17&&\cellcolor{Gray2}1292.5~(0.0937)&\cellcolor{Gray3}1291.5&\cellcolor{Gray3}768.00&&\cellcolor{Gray1}1291.45~(0.0124)&\cellcolor{Gray1}1291.45&\cellcolor{Gray1}90.62\\
Golden17&707.76&\cellcolor{Gray2}707.81~(0.0071)&\cellcolor{Gray1}707.76&\cellcolor{Gray2}937.59&&\cellcolor{Gray3}708.06~(0.0424)&\cellcolor{Gray3}707.79&\cellcolor{Gray3}1770.60&&\cellcolor{Gray1}707.76~(0.0000)&\cellcolor{Gray1}707.76&\cellcolor{Gray1}194.76\\
Golden13&857.19&\cellcolor{Gray2}860.00~(0.3278)&\cellcolor{Gray1}857.19&\cellcolor{Gray2}910.35&&\cellcolor{Gray3}862.95~(0.6720)&\cellcolor{Gray3}861.02&\cellcolor{Gray3}1749.60&&\cellcolor{Gray1}858.34~(0.1342)&\cellcolor{Gray1}857.19&\cellcolor{Gray1}330.77\\
Golden9&579.71&\cellcolor{Gray3}585.21~(0.9488)&\cellcolor{Gray2}583.24&\cellcolor{Gray3}1720.76&&\cellcolor{Gray2}584.88~(0.8918)&\cellcolor{Gray3}583.48&\cellcolor{Gray1}1427.40&&\cellcolor{Gray1}580.70~(0.1708)&\cellcolor{Gray1}580.06&\cellcolor{Gray2}1434.62\\
Golden18&995.13&\cellcolor{Gray3}997.85~(0.2733)&\cellcolor{Gray2}995.65&\cellcolor{Gray2}2297.62&&\cellcolor{Gray2}997.53~(0.2412)&\cellcolor{Gray3}996.1&\cellcolor{Gray3}3447.60&&\cellcolor{Gray1}997.22~(0.2100)&\cellcolor{Gray1}995.14&\cellcolor{Gray1}428.27\\
Golden14&1080.55&\cellcolor{Gray1}1082.15~(0.1481)&\cellcolor{Gray1}1080.55&\cellcolor{Gray2}1513.32&&\cellcolor{Gray3}1084.78~(0.3915)&\cellcolor{Gray3}1082.08&\cellcolor{Gray3}3148.20&&\cellcolor{Gray2}1082.78~(0.2064)&\cellcolor{Gray1}1080.55&\cellcolor{Gray1}393.38\\
Golden10&735.43&\cellcolor{Gray3}744.17~(1.1884)&\cellcolor{Gray3}741.96&\cellcolor{Gray3}3229.35&&\cellcolor{Gray2}743.57~(1.1068)&\cellcolor{Gray2}739.6&\cellcolor{Gray2}2406.60&&\cellcolor{Gray1}739.05~(0.4922)&\cellcolor{Gray1}736.86&\cellcolor{Gray1}1662.97\\
Golden19&1365.6&\cellcolor{Gray2}1367.25~(0.1208)&\cellcolor{Gray2}1366.29&\cellcolor{Gray2}2917.31&&\cellcolor{Gray3}1367.62~(0.1479)&\cellcolor{Gray3}1366.36&\cellcolor{Gray3}5731.20&&\cellcolor{Gray1}1366.33~(0.0535)&\cellcolor{Gray1}1365.60&\cellcolor{Gray1}412.98\\
Golden15&1337.27&\cellcolor{Gray3}1349.23~(0.8944)&\cellcolor{Gray3}1347.13&\cellcolor{Gray2}3265.68&&\cellcolor{Gray2}1347.79~(0.7867)&\cellcolor{Gray2}1343.55&\cellcolor{Gray3}6049.80&&\cellcolor{Gray1}1339.96~(0.2012)&\cellcolor{Gray1}1338.00&\cellcolor{Gray1}1244.74\\
Golden11&911.98&\cellcolor{Gray3}922.93~(1.2007)&\cellcolor{Gray3}921.46&\cellcolor{Gray3}5978.97&&\cellcolor{Gray2}920.72~(0.9584)&\cellcolor{Gray2}917.24&\cellcolor{Gray2}4755.00&&\cellcolor{Gray1}916.10~(0.4518)&\cellcolor{Gray1}913.38&\cellcolor{Gray1}2188.99\\
Golden20&1817.59&\cellcolor{Gray3}1823.52~(0.3263)&\cellcolor{Gray3}1821.16&\cellcolor{Gray2}4997.31&&\cellcolor{Gray2}1820.64~(0.1678)&\cellcolor{Gray2}1820.11&\cellcolor{Gray3}7380.60&&\cellcolor{Gray1}1819.18~(0.0875)&\cellcolor{Gray1}1818.11&\cellcolor{Gray1}1033.22\\
Golden16&1611.28&\cellcolor{Gray3}1627.76~(1.0228)&\cellcolor{Gray3}1624.55&\cellcolor{Gray2}4835.12&&\cellcolor{Gray2}1624.47~(0.8186)&\cellcolor{Gray2}1616.13&\cellcolor{Gray3}9621.00&&\cellcolor{Gray1}1615.19~(0.2427)&\cellcolor{Gray1}1612.44&\cellcolor{Gray1}3005.96\\
Golden12&1101.5&\cellcolor{Gray3}1116.52~(1.3636)&\cellcolor{Gray3}1113.30&\cellcolor{Gray3}10410.70&&\cellcolor{Gray2}1112.26~(0.9768)&\cellcolor{Gray2}1108.78&\cellcolor{Gray2}7251.00&&\cellcolor{Gray1}1106.18~(0.4249)&\cellcolor{Gray1}1102.64&\cellcolor{Gray1}3751.89\\
\hline
\end{tabular}
  \begin{tablenotes}
\item[1] Intels CoreTM i7 with 2.93 GHz and 8 GB of RAM running under Ubuntu Linux 64 bits.
\item[2] Xeon E5-2650 v2 CPU at 2.60 GHz.
\item[3]Intel Xeon E5-2680v2 processor with 2.8 GHz, 10 cores and 128 GB DDR3 1866 MHz RAM.
\end{tablenotes}
\end{threeparttable}
\end{table}

The results indicate that, except for instance Golden14, which it found the second best result, AILS-PR obtained best average solutions (\textit{gaps}) in all tested instances. Moreover, it achieved the top best solutions in all instances. The computational time was also the best in almost all instances, being slower than ILS-SP and SISRs in two and one instances, respectively.

\bibliographystyle{elsarticle-harv} \bibliography{Bibliography}